\documentclass[a4paper, 11pt]{amsart}
\usepackage{amsmath,amssymb,amsthm,amscd,latexsym}
\usepackage[all]{xy}
\usepackage[dvips]{graphicx}
\usepackage{hyperref}

\input xypic
\xyoption{all}

\newcommand{\rar}{\rightarrow}
\newcommand{\lar}{\longrightarrow}
\newcommand{\llar}{-\kern-5pt-\kern-5pt\longrightarrow}

\newtheorem{Theorem}{Theorem}[section]
\newtheorem{Lemma}[Theorem]{Lemma}
\newtheorem{Corollary}[Theorem]{Corollary}
\newtheorem{Proposition}[Theorem]{Proposition}
\newtheorem{Remark}[Theorem]{Remark}
\newtheorem{Example}[Theorem]{Example}
\newtheorem{Conjecture}[Theorem]{Conjecture}
\newtheorem{Definition}[Theorem]{Definition}
\newtheorem{Question}[Theorem]{Question}

\def\sqr#1#2{{\vcenter{\hrule height.#2pt
        \hbox{\vrule width.#2pt height#1pt \kern#1pt
            \vrule width.#2pt}
        \hrule height.#2pt}}}
\def\phi{\varphi}
\def\demo{\noindent{\bf Proof. }}
\def\square{\mathchoice\sqr64\sqr64\sqr{4}3\sqr{3}3}
\def\qed{\hspace*{\fill} $\square$}

\DeclareMathOperator{\depth}{depth}

\DeclareMathOperator{\A}{\alpha}

\DeclareMathOperator{\D}{d}

\def\xx{{\bf x}}

\def\TT{{\bf T}}

\def\bb{\mathbf{b}}

\def\fm{{\mathfrak m}}

\def\Ree#1{{\mathcal R}(#1)}

\def\depth{{\rm depth}\,}

\def\ker{{\rm ker}\,}

\def\restr{{\kern-1pt\restriction\kern-1pt}}

\def\C{{\mathbb C}}

\def\pp{{\mathbb P}}

\begin{document}
\title[The ubiquity of Sylvester forms]{The ubiquity of Sylvester forms in almost complete intersections}

\author{Aron Simis}
\address{Departamento de Matem\'atica \\
Universidade Fed. de Pernambuco \\
50740-540 Recife, PE\\
Brazil}
\address{\rm and}
\address{Departamento de Matem\'atica\\
Universidade Fed. da Paraiba\\
58051-900 J. Pessoa, PB,
Brazil}
\email{aron@dmat.ufpe.br}

\author{Stefan O. Toh\v aneanu}
\address{Department of Mathematics\\
  University of Idaho\\
Moscow, Idaho 83844-1103, USA}
\email{tohaneanu@uidaho.edu}

 \keywords{Rees algebra, Sylvester forms, almost Cohen--Macaulay, reduction number, monomials, birational. \\ \indent The first author is partially supported by a CNPq grant and a PVNS Fellowship from CAPES.}







\begin{abstract}
The subject matter is the structure of the Rees algebra of almost complete intersection ideals of finite colength in low-dimensional polynomial rings over fields.
The main tool is a mix of Sylvester forms and iterative mapping cone construction.
The material developed spins around ideals of forms in two or three variables in the search of those classes for which the corresponding Rees ideal is generated by Sylvester forms and is almost Cohen--Macaulay.
A main offshoot is in the case where the forms are monomials.
Another consequence is a proof that the Rees ideals of the base ideals of certain plane Cremona maps (e.g., de Jonqui\`eres maps) are generated by Sylvester forms and are almost Cohen--Macaulay.

\end{abstract}

\maketitle


\section*{Introduction}

Our goal  is the study of the Rees algebras of classes of almost complete intersection ideals of finite colength in polynomial rings over fields.
The emphasis is on those aspects of the structure of the presentation ideal of the Rees algebra  that may benefit from the appeal to Sylvester forms.
A crucial supporting role will be deployed by a technique of iterated mapping cone construction.

A major focus of this work is on the Rees algebra of a monomial ideal in few variables. Our main source of inspiration for the questions envisaged here comes from the last subsection of \cite{syl3}, particularly Conjecture 4.15.
However, the methods employed in this work are different from those in that reference.

Let  $R=k[x_1, \ldots, x_n]$ denote a polynomial ring over a field $k$.
By an {\em almost complete intersection} we mean an ideal
$I=(a_1, \ldots, a_g, a_{g+1})\subset R$ of codimension $g$ where the subideal $\{a_1, \ldots, a_g\}$ is a
regular sequence generating an ideal not containing  $a_{g+1}$.
Of great interest is the case where $a$'s are forms of the same degree. Except for a few applications to ternary codimension $2$ ideals, $I$ will be an ideal of finite colength.

The Rees algebra $\Ree I=R[IT]\subset R[T]$ has a presentation ideal
$$\mathcal I:=\ker(R[t_1,\ldots,t_{g+1}] \lar R[IT]),\; t_j\mapsto a_jT,$$
often informally referred to as the {\em Rees ideal} of $I$.
 There is a non dismissive importance in asking about the nature of a minimal set of generators of $\mathcal I$ (or at least about the bidegrees of its generating forms in the case $I$ is generated in fixed degree -- see, e.g., \cite{aha}, where this has bearing to the question of birationality).
For example, as is well-known, a substantial subset $L$ of generators come from the syzygies of $I$.
We will focus on several situations where $\mathcal I$ can be generated by those plus a set of Sylvester forms derived from $L$ by iteration.
This sort of goal has been  previously obtained in a few limited situations that will be used here.
To extend this sort of result to other environment is the content of two of the main theorems in this work (Theorem~\ref{main_binary} and Theorem~\ref{main3} (a)).

While the nature of sets of generators of $\mathcal I$ is pertinent, our interest often resides in properties of $\Ree I$ as an $R$-algebra, e.g. its depth and other metrics such as the reduction number of $I$ (finite colength situation).
A major question is a lower bound for the depth.
In this regard one asks what happens beyond the inequality
${\rm depth}(\Ree I)\leq \dim \Ree I=n+1$, where the depth is computed on the maximal graded ideal $(\fm, \Ree I_+)$, with $\fm=(x_1,\ldots,x_n)$.
One says that $\Ree I$ is {\em almost Cohen--Macaulay} (respectively, {\em strictly almost Cohen--Macaulay}) if ${\rm depth}(\Ree I)\geq n$ (respectively, ${\rm depth}(\Ree I)= n$).

In the following situations of almost complete intersections we provide a neat description of a set of generators of $\mathcal I$ and, in addition, show that  $\Ree I$ is  almost Cohen--Macaulay:

$\bullet$ (non $\fm$-primary) $I\subset k[x,y,z]$ is the base ideal of a plane de Jonqui\`eres map

$\bullet$ (finite colength) $I\subset k[x,y]$ is generated by monomials of same degree $x^a,y^a,x^by^{a-b}$, with $1\leq b\leq a$

$\bullet$ (finite colength) $I\subset k[x,y,z]$ is generated by monomials of the form $x^a,y^a,z^a, (xyz)^b$, with $1\leq b <a$.

The part of the almost Cohen--Macaulayness in the second of these results was proved by Rossi--Swanson (\cite[Proposition 1.9]{RoSw}) in more generality, using the Ratliff--Rush filtration. However, their proof scheme gives no insight whatsoever into the structure of $\mathcal I$ whereas ours includes a detailed description of an iterative set of generators of $\mathcal I$ which are binomial Sylvester forms,  in such a manner as to yield at each step a free resolution by a mapping cone construction.
The last step gives a free resolution of $\mathcal I$, from which one sees the almost Cohen--Macaulay fact.

The limitation of the use of the  Ratliff--Rush filtration is that it has no impact on the almost Cohen--Macaulay condition beyond dimension $2$. The device in this work, technical as it resembles, seems to be a lone alternative.
Anyway, it works in the third result above and opens a glade into how to get the case of $n$ variables for this class of monomial ideals.

The binary environment (i.e., $R=k[x,y]$) allows for some simplification in the case of forms of the same degree as the rational map defined by these forms is  birational onto the image provided the forms are not reparameterizable. As a consequence, the value of the first Hilbert coefficient $e_1(I)$ is known. In addition, if the forms are simple enough (like monomials or binomials) then the lengths of the modules $I^r/JI^{r-1}$ ($J$ a reduction) are under some control.
In this situation one can efficiently use the criterion of Huckaba--Marley to decide almost Cohen--Macaulayness.

We spent some time describing these lengths in the case of monomials. The results are unexpectedly rich and some of the questions remain open.

In the case of  $n\geq 3$ variables, even the
the harmless looking case of monomials of the same degree gives a sort of punch line in the theory since the corresponding rational map is not birational onto the image.
In particular there is no a priori formula to get $e_1(I)$ -- actually, granting the almost Cohen--Macaulay property and some flirting around with the lengths of the modules $I^r/JI^{r-1}$ gives a way to compute $e_1(I)$.

\medskip

We now briefly describe the contents of each section, by stressing that throughout the main actor is an almost complete intersection of finite colength.

\smallskip

In the first section we bring up the tools to be used in many parts of the paper. Sylvester forms are introduced in a brushlike way, along with their natural relation to the construction of certain mapping cones.  The main thread is a reworking of a binary ideal that admits a linear syzygy and its consequences to the homological structure of its Rees algebra and of the Rees algebra of the base ideal of a plane (i.e., in $\pp^2$) de Jonqui\`eres map.
The second of these results will come as a consequence of the first through a  specialization procedure from plane Cremona maps to birational maps $\pp^1\dasharrow \pp^2$.
As a corollary, we show that the Rees algebra of the base ideal of a plane de Jonqui\`eres map is almost Cohen--Macaulay, thus completing the homological side of such maps as given in \cite[Theorem 2.7 (iii)]{HS}.

\smallskip

The second section deals with monomial almost complete intersections of finite colength.
We first explore general aspects of this setup, such as the question of reparametrization and minimal reductions with the corresponding reduction numbers.
Obtaining a precise information on the latter seems to incur in an integer optimization problem.
Definite results are obtained in the case that we call {\em uniform} (first treated in \cite{syl3}), namely, when $I=(x_1^{a},\ldots, x_n^{a}, (x_1\cdots x_n)^b)$, with $0<b<a$.

We next consider binary ideals generated by monomials of the same degree. The main result is that, for such an ideal, its Rees ideal $\mathcal I$  is generated by binomial Sylvester forms starting out from the syzygy generators. Coupled with the technique of iterated mapping cones it gives a proof of the almost Cohen--Macaulayness of $\mathcal I$.
The proof of the Sylvester like structure of the latter is unavoidably technical. Of course, getting some examples in the computer gives one a sense of how this structure comes about, but writing a rigorous proof is altogether a different matter.

An appropriate technique to organize the proof spins around the Euclidean algorithm as an iterative procedure.
This tool has been independently used by Benitez--D'Andrea (\cite{Dandrea})
who also found a set of generators of the Rees ideal $\mathcal I$. Whereas our use of it appeals to certain Fibonacci like sequences,  their approach bring up so-called {\em extended Euclidean reminder sequences}.
Since both approaches involve a good deal of calculations, we have decided to stick to our own as, being more familiar to us, leads straight to the present goals.
In any case, the priority goes to their approach having been first to come public.

We have complemented these results giving the details of the Huckaba--Marley environment,  calculating the lengths of the various modules $I^r/JI^{r-1}$, with $J\subset I$ a reduction.
Since the Hilbert coefficient $e_1(I)$ is known a priori in this situation, in principle the almost Cohen--Macaulay property can be verified directly. Beyond this, there are some structural questions left open about the linear syzygies of the powers of $I$.

The last part of the section  deals with uniform case in three variables.
Note that the uniform binary case is devoided of interest since, up to reparametrization, it reduces to the ideal $(x^2,y^2,xy)$, but not so the ternary case.
Resorting to similar techniques as before, we show that the Rees ideal $\mathcal I$ is generated by the syzygy forms and binomials Sylvester forms, and that the Rees algebra is almost Cohen--Macaulay.

\section{General methods}

In this section we focus on the structure of the Rees ideal of  in that it benefits from the two techniques mentioned in the Introduction: Sylvester forms and mapping cones.

\subsection{Sylvester forms}

We recall the notion of a Sylvester form in the special case to be used in this work (for additional details see \cite[Section 2]{syl1}, \cite[Section 4]{syl2}.

\begin{Definition}\rm
Let $(\alpha,\beta)\subset R=k[x_1,\ldots,x_n]$ be an ideal generated by two nonzero forms and let $f,g\in (\alpha,\beta)R[t_1,\ldots,t_m]$ be given biforms. Write this containment as a matrix equality
$$\left[\begin{array}{c} f \\ g \end{array}\right]=
                     \left[\begin{array}{ll} c_1 & c_2 \\
                      d_1 & d_2 \end{array}\right]
                      \left[\begin{array}{c}\alpha
                      \\ \beta \end{array} \right],
$$
for suitable biforms $c_i,d_i\in  S=R[t_1,\ldots,t_m]$.

The {\em Sylvester form} of $f,g$ with respect to $(\alpha,\beta)$ is the determinant of the above content $2\times 2$ matrix,
 denoted $\det(f,g)_{(\alpha,\;\beta)}$.
\end{Definition}

Our main use of Sylvester forms is in the case where $f,g$ are biforms in the Rees  ideal $\mathcal I\subset S$.
Under this assumption, by Cramer one has
$$\det(f,g)_{(\alpha,\;\beta)}\cdot (\alpha,\beta)\subset (f,g)\subset \mathcal I.$$
Since $\mathcal I$ is a prime ideal and $\mathcal I\cap R=\{0\}$, the determinant belongs to $\mathcal I$.

Our first example of the joint use of these techniques is a different proof of a basic result of algebraic nature in elimination theory,
carrying additional information in the homological side.

 In a slightly different form, but in a more encompassing environment, the following result can be found in \cite[Corollary A.140]{Vas} and goes back to Northcott.
 It will be applied in the special environment of Sylvester forms.

\begin{Lemma}\label{mapcone} Let $S$ be a  commutative ring  and let  $\{A,B\}$ and $\{C,D\}$ be two regular sequences. Let $a,b,c,d\in S$, be given such that  $$\left[\begin{array}{c}C\\D\end{array}\right]=\underbrace{\left(\begin{array}{cc}a & b\\c& d \end{array}\right)}_{\mathcal M}\left[\begin{array}{c}A\\B\end{array}\right].$$
One has:
\begin{enumerate}
\item If some entry of $\mathcal M$ is a nonzero divisor modulo $(A,B)$, then $(C,D):E=(A,B),$
where $E:=\det(\mathcal M)$.
\item The mapping cone given by the map of complexes $$\begin{array}{lllllllll}
0 & \rar & S & \stackrel{\tiny{\left[\begin{array}{r}-D\\C\end{array}\right]}}\lar & S^2 & \stackrel{\tiny{\left[\begin{array}{cc}C&D\end{array}\right]}}\lar & S &  \rar & 0\\
&& \uparrow\small{\shortparallel} && \uparrow \small{\mathcal M^{*}} && \uparrow \small{\cdot E} && \\
0 & \rar & S & \stackrel{\tiny{\left[\begin{array}{r}-B\\A\end{array}\right]}}\lar & S^2 & \stackrel{\tiny{\left[\begin{array}{cc}A&B\end{array}\right]}}\lar & S &  \rar & 0
\end{array},
$$
is a free resolution of $S/(C,D,E)$, where $\mathcal M^{*}:=\left(\begin{array}{rr}d & -c\\-b& a \end{array}\right)$.
\end{enumerate}
\end{Lemma}
\demo Since $\left[\begin{array}{cc}C&D\end{array}\right]=\left[\begin{array}{cc}A&B\end{array}\right]\cdot\mathcal M^t$, then $$\left[\begin{array}{cc}C&D\end{array}\right]\cdot \mathcal M^{*}=E\cdot \left[\begin{array}{cc}A&B\end{array}\right],$$ and therefore $(A,B)\subseteq (C,D):E$.

Let $f\in (C,D):E$ and, say, $c$ is a non-zero divisor in $S/(A,B)$.

There exist $\alpha,\beta\in S$ such that $$\underbrace{(ad-bc)}_{E}f=\alpha C+\beta D.$$
Multiplying this equation by $B$ and from the fact that $C=aA+bB, D=cA+dB$, we have $$\left(a(D-cA)-c(C-aA)\right)\,f=\alpha BC+\beta BD.$$
This leads to $$D(af-\beta B)=C(\alpha B+cf).$$ Since $\{C,D\}$ is a regular sequence, there exists $\gamma\in S$ such that $\alpha B+cf=\gamma D$. Since $D\in (A,B)$, we have $$f\in(A,B):c=(A,B),$$ as $c$ is a non-zero divisor modulo $(A,B)$.

Therefore $(A,B)=(C,D):E$. From this, the short exact sequence $$\begin{array}{ccccccccc}
0 & \rar & \frac{S}{(C,D):E} & \stackrel{\tiny{\cdot E}}\lar & \frac{S}{(C,D)} & \lar & \frac{S}{(C,D,E)} &  \rar & 0\end{array}$$ shows in the usual way that the mapping cone is free resolution of
the rightmost module.
\qed

\medskip

One knows that $S/(C,D,E)$ is a perfect module of codimension $2$, hence the mapping cone above gives a non-minimal free resolution in this case. However, the result will express the details of the map between the two complexes which are used in the proof of the next proposition.

\medskip

\begin{Proposition}\label{resolution}{\rm (\cite[Theorem 4.4, Proposition 4.6]{syl3})}
Let $I\subset R=k[x,y]$ denote an almost complete intersection of finite colength, generated in fixed degree $d\geq 2$. If $I$ admits a linear syzygy then the Rees algebra $\mathcal{R}(I)$ is almost Cohen--Macaulay.
Moreover, it is strictly so if and only if $d\geq 3$, in which case its minimal free resolution over
the polynomial ring $B=R[t,u,v]$ has the form
\begin{equation}\label{free_resolution}
0\rar B^{d-2}\lar B^{2(d-1)}\lar B^{d+1} \lar B\lar \mathcal{R}(I)\rar 0
\end{equation}
\end{Proposition}
\demo The proof is by induction on the number $d+1$ of minimal generators of the presentation ideal
$\mathcal{I}$ in a presentation of $\mathcal{R}(I)$ over $B$.
We draw upon the result that
\begin{equation}\label{Rees_eqs_binary}
\mathcal{I}=(F,G,H_1, H_2 \ldots, H_{d-2},E),
\end{equation}
where these generators are of respective bidegrees $(1,1), (d-1,1), (d-2,2), \ldots, (1,d-1), (0,d)$ (see \cite{CHW}, \cite{syl3}).
Moreover, we will use that these generators are Sylvester forms obtained as determinants
of the successive contents relative to the maximal ideal $(x,y)$.
In particular, $(x,y)H_1\subset (F,G)$ and, in general, $(x,y) H_i\subset (F,G,H_1, H_2 \ldots, H_{i-1})$
for $i=2,\ldots, d-2$ and, finally, $(x,y) E\subset (F,G,H_1, H_2 \ldots, H_{d-2})$.

For $d=2$ we use Lemma~\ref{mapcone} with $A=F,B=G$ and $\{C,D\}=\{y,x\}$, whereby the mapping cone gives a free resolution
of  $\mathcal{I}=B/(F,G,H_1)$ and the latter is
 Cohen--Macaulay.

Suppose $d\geq 3$.
Consider the map of complexes induced by the multiplication map $B\rar B$ by $H_2$
$$\begin{array}{ccccccccc}
0 & \rar & B^2 & \lar & B^3 & \lar & B &  \rar & 0\\
&& \uparrow && \uparrow && \uparrow && \\
0 & \rar & B & \lar & B^2 & \lar & B &  \rar & 0
\end{array},
$$
where the upper complex is the resolution of $B/(F,G,H_1)$
and the lower one is the Koszul complex of $\{y,x\}$.

The resulting mapping cone
\begin{equation}\label{degree4}
0\rar B\lar B^4\lar B^4\lar B
\end{equation}
is a free resolution of $B/(F,G,H_1, H_2)$.
Suppose shown that the free resolution of
$$B/(F,G,H_1, \ldots, H_{i-1})$$
has the form
$$0\rar B^{i-2}\lar B^{2(i-1)}\lar B^{i+1} \lar B.
$$
Then the mapping cone of the map of complexes
$$\begin{array}{ccccccccccc}
0 & \rar & B^{i-1} & \lar & B^{2i} & \lar & B^{i+2} & \lar & B & \rar & 0\\
&&\uparrow && \uparrow && \uparrow && \uparrow && \\
&&0 & \rar & B & \lar & B^2 & \lar & B &  \rar & 0
\end{array},
$$
induced by multiplication $B\rar B$ by $H_i$ is a free resolution
$$0\rar B^{i-1}\lar B^{2i}\lar B^{i+2} \lar B
$$
of $B/(F,G,H_1, \ldots, H_{i})$.
Proceeding this way we arrive at the last map of complexes
$$\begin{array}{ccccccccccc}
0 & \rar & B^{d-3} & \lar & B^{2(d-2)} & \lar & B^{d} & \lar & B & \rar & 0\\
&&\uparrow && \uparrow && \uparrow && \uparrow && \\
&&0 & \rar & B & \lar & B^2 & \lar & B &  \rar & 0
\end{array},
$$
where the top complex is a resolution of $B/(F,G,H_1, \ldots, H_{d-2})$
and the rightmost vertical map is multiplication by the implicit equation $E$.
The mapping cone is a resolution of the entire Rees algebra.
\qed

\subsection{Specialization to the binary case}

We need one additional result related to the mapping cone construction.

\begin{Lemma}\label{depth_equality}
Let $I\subset R:=k[x,y]$ denote an almost complete intersection of finite colength, generated in fixed degree $d\geq 3$.
Let $\mathcal{I}$ stand for the defining ideal of the Rees algebra $\mathcal{R}_R(I)$ of $I$ on the polynomial ring $B:=R[t,u,v]$, and let $\mathcal{K}\subset \mathcal{I}$ denote the subideal contained in $(x,y)B$.
Assume that the rational map defined by the generators of $I$ is birational onto the image.
Then
\begin{equation}
{\rm depth}\,  \mathcal{R}_R(I)=
\left\{
\begin{array}{ll}
{\rm depth} \,B/\mathcal{K}, & \mbox{\rm if} \;\; {\rm depth} \,B/\mathcal{K}\leq 2\\
2, & \mbox{\rm if}\;\; {\rm depth} \,B/\mathcal{K}=3
\end{array}
\right.
\end{equation}
\end{Lemma}

\demo
By assumption, one has $\mathcal{R}_R(I)\simeq B/(\mathcal{K},E)$, where $E$ denotes the implicit equation.
It also clear, by the definition of $\mathcal{K}$, that $E\, (x,y)B\subset \mathcal{K}$.

Let
$$0 \rar B^{b_r} \lar\cdots\lar  B^{b_3} \lar  B^{b_2} \lar B^{b_1} \lar B$$
stand for the minimal graded resolution of $B/\mathcal{K}$ over the standard graded polynomial ring $B=k[x,y,t,u,v]$, where $2\leq r\leq 5$.

Tensor the minimal $R$-resolution of $R/(x,y)$ with $B$ over $R$ and consider the map of complexes induced by multiplication by $E$ on $B$:
$$\begin{array}{ccccccccccccc}
0 & \rar & B^{b_r} & \lar & \cdots & \lar & B^{b_3} & \lar & B^{b_2} & \lar & B^{b_1} & \lar & B\\
&&&&&& \uparrow & & \uparrow & & \uparrow & &  \quad\uparrow \cdot E\\ 
&&&&&& 0 & \lar & B & \lar & B^2 & \stackrel{(x\; y)}{\lar} & B
\end{array}.$$
There are two possibilities according to whether $r\geq 3$ or $r=2$.
If $r\geq 3$ the resulting mapping cone gives the following resolution of $B/(\mathcal{K},E)= \mathcal{R}_R(I)$:
$$0 \rar B^{b_r} \lar\cdots\lar  B^{b_4}\lar  B^{b_3+1} \lar  B^{b_2+2} \lar B^{b_1+1} \lar B.$$
Therefore, the projective dimension of $\mathcal{R}_R(I)$  is at most $r$, hence its depth is at least ${\rm depth} \,B/\mathcal{K}$.

If $r=2$,  i.e., $B/\mathcal{K}$ is Cohen--Macaulay, then the result is the  resolution
$$0\rar B\lar B^{b_2+2} \lar B^{b_1+1} \lar B$$
of $ \mathcal{R}_R(I)$.
In this case,  ${\rm depth}\mathcal{R}_R(I)\leq {\rm depth} \,B/\mathcal{K}$.

The first case gives strict equality of depths  provided the corresponding resolution is minimal; similarly, in the case of $r=2$ if the corresponding resolution is minimal then the depth of $\mathcal{R}_R(I)$ is exactly one less.

To avoid such a direct verification, one may proceed as follows.

{\bf Claim:} $\mathcal{K}\cap EB= E\, (x,y)B$.

The inclusion $E\, (x,y)B\subset (\mathcal{K}\cap EB$ is clear.
For the inclusion $\mathcal{K}\cap EB\subset E\cdot (x,y)$ let $G\in B$ be such that
$GE\in\mathcal{K}\subset (x,y)B$; then $G$ has no term in $t,u,v$ with coefficient in $k$ since
evaluating $x\mapsto 0, y\mapsto 0$ would say that $E$ is a zero divisor in $B$.
This proves the claim.

{\bf Claim:} Multiplication by $E$ on $B$ induces an exact sequence

\begin{equation}\label{easy_sequence}
0\lar B/(x,y)B\stackrel{\cdot E}{\lar} B/\mathcal{K}\lar B/(\mathcal{K}, E)\rar 0.
\end{equation}

Namely, the cokernel is clear, while the kernel is

$$\left(\mathcal{K}\cap  EB, (x,y)B\right)/(x,y)B=\{0\}$$
by the previous claim.

Now, in the case $r\geq 3$, from (\ref{easy_sequence}) and since the left most module is Cohen--Macaulay, the depth of $B/\mathcal{K}$ cannot be strictly less than the depth of $B/(\mathcal{K}, E)$ as otherwise the latter would be  Cohen--Macaulay which is not the case as $r>2$.
Therefore, this gives equality in this case.

Now suppose $r=2$.
In this case we have to show that $\mathcal{R}_R(I)$ is not Cohen--Macaulay.
We can assume that $I=(J,f)$, where $J$ is a reduction of $I$.
If $\mathcal{R}_R(I)$ is Cohen--Macaulay then $I^2=JI$ (see, e.g., \cite[Theorem 2.1]{syl1}).
By a well-known result, the length of $I/J$ is the first Hilbert coefficient of $I$.
Since the corresponding map is assumed to be birational, the latter is ${d \choose 2}$.
On the other hand, $I/J\simeq (J,f)/J\simeq R/J:f$.
Since $J:f$ is generated by two elements of degree $d\geq 3$, the length of $R/J:f$ exceeds ${d \choose 2}$.
This is a contradiction.

This completes the proof.
\qed

\bigskip

We apply this result, via a specialization, to the base ideal of a
Cohen--Macaulay plane Cremona map.
The following theorem does not seem to appear in the related literature in this explicit form.

\begin{Theorem}\label{3vars2vars}
Let $J\subset S:k[x,y,z]$ denote the base ideal of a  Cremona map of  $\pp^2$ of degree $d\geq 2$.
Assume that $S/J$ is Cohen--Macaulay of dimension $1$ and $\ell\in S$ is a $1$-form which is a nonzerodivisor on $S/J$. Then letting $R:=S/(\ell)$ and $I:=(J,\ell)/(\ell)\subset R$, the following hold:
\begin{enumerate}
\item[{\rm (i)}] $I$ is the base ideal of a birational map of $\pp^1$ onto a plane curve of degree $d$.
\item[{\rm (ii)}]
 ${\rm depth} \, \mathcal{R}_S(J)\geq {\rm depth} \, \mathcal{R}_R(I)+1\,${\rm ;} moreover, $\mathcal{R}_S(J)$ is almost Cohen--Macaulay if and only if $\mathcal{R}_R(I)$ is almost Cohen--Macaulay.
\end{enumerate}
\end{Theorem}
\demo
Up to a projective change of coordinates, we may assume that $\ell=z$, in other words, the base locus of the map has no points on the line at infinity.

(i) Since $z$ is regular modulo $J$ and evaluation commutes with taking determinants, we know that $I$ is also the ideal of $2$-minors of the $3\times 2$ matrix over $R$ obtained from the structural syzygy matrix of $J$ by setting $z\mapsto 0$ in all entries.
Thus, $I\subset R$ is an ideal of finite colength generated by $3$ form of degree $d$.

Notice, moreover, that the $k$-subalgebra $k[I_d]\subset R$ has dimension $2$ -- this is because, since $k$ is infinite and $I$ is equi-generated then one can assume that the linear system $I_d$ has $2$ elements forming a regular sequence, and hence are algebraically independent over $k$.

On the other hand,
since $\mathcal{R}_S(J)$ is a domain, $z$ is also regular thereof.
Then, since polynomial relations commute with evaluation, one gets an inclusion
$$(\mathcal{J},z)/(z)  \subset \mathcal{I},$$
up to identification $S[t,u,v]/(z)= R[t,u,v]$, where $\mathcal{J}$ (respectively, $\mathcal{I}$) denote the defining ideal of the Rees algebra $\mathcal{R}_S(J)$ (respectively, $\mathcal{R}_R(I)$).

Now, since $J$ is the base ideal of a birational map of $\pp^2$ then
$\mathcal{J}$ admits enough independent forms of bidegree $(1,s)$ on $k[x,y,z,t,u,v]$, for some $s\geq 1$ (see \cite{aha} for a sufficient explanation of this phenomenon). By setting $z\mapsto 0$ yields enough independent forms of bidegree $(1,s)$ on $k[x,y,t,u,v]$. Since $\dim k[I_d]=2$, this shows (again by \cite{aha}) that the map given by the generating forms of $I$ defines a birational map onto the image; in particular, the latter is a curve of degree $d$.

\medskip

(ii) Let ${\mathcal L}_J$ (respectively, ${\mathcal L}_I$)
denote the defining ideal of the symmetric algebra of $J$ (respectively, $I$).
Then the commuting property of evaluating $z\mapsto 0$ explained in the proof of (i) implies  that $({\mathcal L}_J,z)/(z)\simeq {\mathcal L}_I$ or, equivalently, $({\mathcal L}_J,z)=({\mathcal L}_I,z)$.

By the same token, write $(\mathcal{J}, z)=(\mathcal{A}, z)$, where $\mathcal{A}\subset B:=R[t,u,v]=k[x,y,t,u,v]$ is the uniquely determined subset not involving $z$, i.e., whose elements are the forms obtained by setting $z\mapsto 0$ in every form of a minimal set of homogeneous generators of $\mathcal{J}$.
Set  $\mathcal{A}B$ for the ideal of $B$ it generates.
Clearly, one has inclusions ${\mathcal L}_I\subset \mathcal{A}B\subset \mathcal{K}$, where as in the statement of Lemma~\ref{depth_equality} $\mathcal{K}\subset B$ denote the subideal of $\mathcal{I}$ of all forms with coefficients in $(x,y)$.

\smallskip

{\bf Claim:} $\mathcal{A}B = \mathcal{K}$.

\smallskip

To show this, note the following thread of inclusions
$$\mathcal I={\mathcal L}_I:I^{\infty}={\mathcal L}_I:(x,y)^{\infty}\subset \mathcal{A}B:(x,y)^{\infty}\subset {\mathcal K}:(x,y)^{\infty}\subset \mathcal I.$$
Therefore, we get equalities throughout.
Moreover, as $\mathcal I=(\mathcal K, E)$, where $E$ is the equation of the image, it follows that $(\mathcal A B, E)=(\mathcal K, E)$.

This implies that any generator of $\mathcal K$ can be written in terms of those of $\mathcal A B$ and $E$.
But any minimal generator of $\mathcal K$ has bidegree $(r,s)$, with $r+s\leq d=\deg(E)$ and $r\geq 1$, hence necessarily $s<d$, thus implying that $s\leq d-1$. This says that such a form cannot involve effectively $E$. This forces $\mathcal K\subset \mathcal A B$, as was claimed.

\smallskip

Now, the assertion of this item is equivalent to
$\depth\,\mathcal{R}_S(J)/z\,\mathcal{R}_S(J)=\depth \mathcal{R}_R(I)$.

Note that now
$$\mathcal{R}_S(J)/z\,\mathcal{R}_S(J)\simeq
B/\mathcal{A}B= B/\mathcal{K}$$
and that
 $\depth B/\mathcal{K}\geq \depth \mathcal{R}_R(I)$ by Lemma~\ref{depth_equality} (note that the statement of (ii) is more or less trivial for $d=2$, as both algebras are Cohen--Macaulay; therefore, we may assume that $d\geq 3$ in order to apply the lemma).

For the last additional statement of item (ii), one implication is immediate. For the reverse implication, if $\mathcal{R}_S(J)$ is almost Cohen--Macaulay then so is $B/\mathcal{K}$ by the above.
Then again from  Lemma~\ref{depth_equality}, the depth of $\mathcal{R}_R(I)$ is $2$.
\qed

\bigskip

As a consequence, in the case of a particular, but very fundamental, Cremona map, one has a definite result.

\begin{Corollary}
Set $S:=k[x,y,z]$ and let $\mathfrak{J}\subset S$ denote the base ideal of a plane de Jonqui\`eres map  of degree $d\geq 2$.
Then, for $d\geq 4$, the  Rees algebra of $\mathfrak{J}$ is a strict almost Cohen--Macaulay ring.
Moreover, it has a minimal free resolution of the form
\begin{equation}\label{free_resolution_Jonq}
0\rar A^{d-3}\lar A^{2(d-2)}\lar A^{d} \lar A,
\end{equation}
where $A=S[t,u,v]=k[x,y,z,t,u,v]$.
\end{Corollary}
\demo
It is well known that the base ideal of a plane de Jonqui\`eres map of degree $d\geq 2$
is Cohen--Macaulay, with structural matrix of the form
$$\phi=
\left(
\begin{array}{cc}
x&p_2\\
-y& p_1\\
0 & q
\end{array}
\right)
$$
up to change of coordinates and for suitable forms $p_1,p_2,q$ -- see, e.g., \cite[The proof of Proposition 2.3]{HS} (actually, in this format $z$ is automatically a nonzerodivisor modulo $J$ since $(x,y)$ is the only associated prime).

Therefore, the assertion on the almost Cohen--Macaulayness follows from Theorem~\ref{3vars2vars} and from Proposition~\ref{resolution} (or its predecessors such as \cite[Theorem 4.4, Proposition 4.6]{syl3}).

For the complementary assertion, by \cite[Theorem 2.7 (iii)]{HS}, $\mathcal{R}_S(J)$ is Cohen--Macaulay if and only if $d\leq 3$.
The shape of the minimal resolution for $d\geq 4$ follows from the mapping cone construction of Lemma~\ref{depth_equality} and from the resolution in Proposition~\ref{resolution}.
\qed


\section{Almost complete intersection monomial ideals}

In the previous sections we focused on arbitrary binary almost complete intersections and its bearing to some plane rational maps of $\pp^2$.
Here we wish to look at some questions in the case the forms are monomials.

Quite generally, if $I$ is an ideal generated by forms of the same degree
then these forms span a linear system defining a rational map
of projective spaces.
A major question is as to whether or when this map is birational onto its image. Knowing this a priori facilitates by quite a bit the theory as was reported in \cite{syl1, syl2, syl3}.

Now, a preliminary obstruction to birationality is a more elementary notion related to the phenomenon that the spanning forms might be expressible in terms of forms of lesser degree.
This is of course a known phenomenon, known as {\em re-parametrization}.

We deploy the algebra behind it in the simpler case where the forms are monomials.

\begin{Lemma}\label{reparametrization}
Let $I=(x_1^{a_1},\ldots ,x_n^{a_n}\,,\,x_1^{b_1}\cdots x_n^{b_n})$
with $0\leq b_i< a_i$ for every $i$ and such that there are at least
two different indices $i,j$ for which $b_i\neq 0, b_j\neq 0$.
Set $d_i=\gcd (a_i,b_i)$ for every $i$ such that $b_i\neq 0$.
Consider the ideal
$$I':=(x_1^{a'_1},\ldots ,x_n^{a'_n}\,,\,x_1^{b'_1}\cdots x_n^{b'_n}),$$
where $a'_i=a_i/d_i, b'_i=b_i/d_i$ if $b_i\neq 0$, and $a'_i=a_i$ otherwise.
Consider the ring endomorphism $\delta$ of $R[\TT]=R[T_1,\ldots,T_{n+1}]$
that sends $T_j\mapsto T_j$ for every $j=1,\ldots,n+1$ and sends $x_i\mapsto x_i^{c_i}$
where $c_i=d_i$ if $b_i\neq 0$, $c_i=1$ otherwise.
Then:
\begin{enumerate}
\item
If $\mathcal{J}'$ is the defining ideal of $\mathcal{R}_R(I')$
then $\delta(\mathcal{J}')$ is the defining ideal of $\mathcal{R}_R(I)$.
Moreover, the two defining ideals have the same number of minimal
binomial generators with the same $\TT$-degrees throughout.
\item If $I'$ has a monomial reduction then so does $I$ and the corresponding
reduction numbers coincide.
\end{enumerate}
\end{Lemma}
\demo
(1) Set $S:=R[\TT]$.
 Let $\mathcal{K}(I):=I_1({\bf T}\cdot \mathcal{K}_2(I))$ (respectively,  $\mathcal{K}(I'):=I_1({\bf T}\cdot \mathcal{K}_2(I'))$) stand for the ideal generated by the set of binomials coming from the Koszul relations of $I$ (respectively,  $I'$), where  $ \mathcal{K}_2(I)$ ((respectively, $\mathcal{K}_2(I'))$) denotes the matrix of the Koszul syzygies of the generators of $I$ (respectively, $I'$).

Some $b_i\neq 0$, so assume that this is the case for $i=1$.
Since $\mathcal{J} = \mathcal{K}(I) : I^{\infty}$ (respectively,
$\mathcal{J}' =  \mathcal{K}(I') : I'^{\infty}$) then it is easy to see that $\mathcal{J} = \mathcal{K}(I) : (x_1^{a_1})^{\infty}$
(respectively, $\mathcal{J}' =  \mathcal{K}(I') : (x_1^{a'_1})^{\infty}$.
On the other hand, it is clear, by definition, that  $\delta(\mathcal{K}(I'))= \mathcal{K}(I)$.
Therefore, $\mathcal{J}  =  \delta( \mathcal{K}(I')) :  (x_1^{a_1})^{\infty}$.

Clearly, then $\delta(\mathcal J')=\delta( \mathcal{K}(I') : (x_1^{a'_1})^{\infty})\subset \mathcal J$.
Conversely, we claim that every binomial $H\in \mathcal J$ is of the form $H=\delta(H')$, for some binomial $H'\in \mathcal J'$.
Thus, let $H\in\mathcal J$. Without loss of generality one can assume $$H=x_1^{\alpha_1}\cdots x_p^{\alpha_p}T_{n+1}^{\gamma} T_{p+1}^{\beta_{p+1}}\cdots T_n^{\beta_n}-x_{p+1}^{\alpha_{p+1}}\cdots x_n^{\alpha_n}T_1^{\beta_1}\cdots T_p^{\beta_p},$$ where $\alpha_i,\beta_j,\gamma$ are non-negative integers (some possibly zero) satisfying the exponent equations:

\begin{eqnarray}
\alpha_s+\gamma b_s&=&\beta_s a_s,\mbox{ for }1\leq s\leq p\nonumber\\
\alpha_t&=&\beta_t a_t+\gamma b_t,\mbox{ for }p+1\leq t\leq n.\nonumber
\end{eqnarray}

If $d_i=\gcd(a_i,b_i)$, then $d_i$ divides $\alpha_i$, for all $1\leq i\leq n$. So $\alpha_i=\alpha_i'd_i$. Simplifying by $d_i$ all of the above exponent equations, one obtains

\begin{eqnarray}
\alpha_s'+\gamma b_s'&=&\beta_s a_s',\mbox{ for }1\leq s\leq p\nonumber\\
\alpha_t'&=&\beta_t a_t'+\gamma b_t',\mbox{ for }p+1\leq t\leq n,\nonumber
\end{eqnarray} meaning that $$H=\delta(H'),$$ where $H'=x_1^{\alpha_1'}\cdots x_p^{\alpha_p'}T_{n+1}^{\gamma} T_{p+1}^{\beta_{p+1}}\cdots T_n^{\beta_n}-x_{p+1}^{\alpha_{p+1}'}\cdots x_n^{\alpha_n'}T_1^{\beta_1}\cdots T_p^{\beta_p}\in \mathcal J'$.

\medskip

To see the assertion about minimal generators, we proceed as follows.
Let $G\in \mathcal J'$ be a binomial minimal generator.
Write $$\delta(G)=g_1\delta(G_1)+\cdots+g_m\delta(G_m),$$ for $G_u$ binomials in $\mathcal J'$ and $g_u\in S$.

Keeping in mind that  $\delta(G)$ is a binomial, let ${\rm in}(\delta(G))=cx_1^{\alpha_1}\cdots x_p^{\alpha_p}T_{n+1}^{\gamma} T_{p+1}^{\beta_{p+1}}\cdots T_n^{\beta_n}$ stand for the leading monomial of $\delta(G)$ under some monomial order on $S$, where $c\in k-\{0\}$. Then this monomial must be multiple of some monomial $M$ in $S$ and another monomial $N$ in $S$ effectively appears in one of the binomials $G_u$. Since $N$ divides ${\rm in}(\delta(G))$, then $$N=x_1^{\bar{\alpha}_1}\cdots x_p^{\bar{\alpha}_p}T_{n+1}^{\bar{\gamma}} T_{p+1}^{\bar{\beta}_{p+1}}\cdots T_n^{\bar{\beta}_n},$$
where $\alpha_s\geq\bar{\alpha}_s\geq 0, \gamma\geq\bar{\gamma}\geq 0$, and $\beta_t\geq\bar{\beta}_t\geq 0$.
Then $$M=cx_1^{\alpha_1-\bar{\alpha}_1}\cdots x_p^{\alpha_p-\bar{\alpha}_p}T_{n+1}^{\gamma-\bar{\gamma}} T_{p+1}^{\beta_{p+1}-\bar{\beta}_{p+1}}\cdots T_n^{\beta_n-\bar{\beta}_n}.$$

From the exponent equations for $\alpha$'s and $\bar{\alpha}$'s, for $1\leq s\leq p$ one has $\alpha_s+\gamma b_s = \beta_s a_s$ and $\bar{\alpha}_s+\bar{\gamma} b_s=\bar{\beta}_s a_s$, giving $$\alpha_s-\bar{\alpha}_s=(\beta_s-\bar{\beta}_s) a_s-(\gamma-\bar{\gamma}) b_s, 1\leq s\leq p.$$ Since $d_s=\gcd(a_s,b_s)$, then $\alpha_s-\bar{\alpha}_s=d_s\rho_s$. So we obtain $$M=\delta(M'),$$ where $M'=cx_1^{\rho_1}\cdots x_p^{\rho_p}T_{n+1}^{\gamma-\bar{\gamma}} T_{p+1}^{\beta_{p+1}-\bar{\beta}_{p+1}}\cdots T_n^{\beta_n-\bar{\beta}_n}$.

Suppose $u=1$. Then, since $\delta$ is a ring endomorphism we got $$\delta(G-M'G_1)=(g_1-M)\delta(G_1)+g_2\delta(G_2)+\cdots+g_m\delta(G_m),$$ with ${\rm in}(\delta(G-M'G_1))<{\rm in}(\delta(G))$. Inducting, one gets $G\in(G_1,\ldots,G_m)$, contradicting the minimality of $G$. So $\delta(G)$ is also a part of a minimal generating set for $\mathcal J$.

(2) Let us assume $J':=(x_1^{a'_1},\ldots, x_n^{a'_n})$ is a reduction for $I'$. Since both $I'$ and $I$ are aci's then the powers behave as: if $I' = (J',f')$ then $I'^r = (J'I'^{r-1}, f'^r)$ , where $f'=x_1^{b_1'}\cdots x_n^{b_n'}$ and if we denote with $J:=(x_1{a_1},\ldots,x_n^{a_n})$, we have $I = (J,f)$ then $I^p = (JI^{p-1}, f^p)$ , where $f=x_1^{b_1}\cdots x_n^{b_n}$.

Suppose ${\rm red}_{J'}(I')=\ell'$ and ${\rm red}_J(I)=\ell$. Then $\ell'$ and $\ell$, respectively, are the smalles integers such that $f'^{\ell'+1}\in J'I'^{\ell'}$ and $f^{\ell+1}\in JI^{\ell}$.

Since everything is monomial, $\delta(f')=f$, $\delta(J')=J$, $\delta(I')=I$, and $\delta$ is a ring endomorphism, one has indeed that $\ell'=\ell$.
\qed

\bigskip

We wish to address the following question  raised in \cite[Section 4.2]{syl3}:
\begin{Question}\rm (\cite[Conjecture 4.15]{syl3})\label{main_question}
Let $I\subset R:=k[x_1, \ldots, x_n]$ denote a monomial ideal which is an almost complete intersection
of finite colength. Is its Rees algebra $\mathcal{R}_R(I)$ almost Cohen--Macaulay?
\end{Question}

We assume throughout that $I$ is not a complete intersection.
Therefore, $I$ is minimally generated by forms $x_1^{a_1},\ldots, x_n^{a_n}, \xx^{\bb}$,
where $\xx^{\bb}=x_1^{b_1}\cdots x_n^{b_n}$, with $0\leq b_i< a_i$ for every $i$ and there are at least
two different indices $i,j$
for which $b_i\neq 0, b_j\neq 0$.

A much less obvious condition is that the subideal $J:=(x_1^{a_1},\ldots, x_n^{a_n})$ be a minimal reduction of $I$.
Such an assumption implies the inequality
$\sum_{i=1}^n \frac{b_i}{a_i} \geq 1$, while otherwise a typical minimal reduction may fail to be
generated by monomials, requiring some binomial generators.

\begin{Lemma}\label{red_no}
Suppose that $J:=(x_1^{a_1},\ldots, x_n^{a_n})$ is a minimal reduction of $I$.
Then the reduction number ${\rm red}_J(I)$ is the least integer $d\geq 1$ such that there exist $t\geq 2$
 distinct indices $i_1,\ldots,i_t\in\{1,\ldots,n\}$ and corresponding positive integers $s_{i_1},\ldots,s_{i_t}$
 with $s_{i_1}+ \cdots +s_{i_t}=d+1$ satisfying the inequalities $(d+1)b_{i_{\ell}}\geq s_{i_{\ell}}a_{i_{\ell}}$
 for $\ell=1,\ldots,t$.
\end{Lemma}
\demo
Since $I=(J,\xx^{\bb})$, then for any $r\geq 1$, one has
$$I^{r+1}=(JI^r,\xx^{(r+1)\bb})=(J^{r+1}, J^r\xx^{\bb},\ldots, J\xx^{r\bb}, \xx^{(r+1)\bb}).$$
Then ${\rm red}_J(I)$ will be the least $r$ such that $\xx^{(r+1)\bb}\in JI^r$.
But note that all the generator blocks of $JI^r$ are monomials, therefore $\xx^{(r+1)\bb}\in JI^r$ if and only if
$\xx^{(r+1)\bb}\in J^{r+1-s}\xx^{s\bb}$ for some $s\in \{0,\ldots,r\}$.
However,  this inclusion is only possible if $s=0$ since otherwise we could cancel a copy of $\xx^{\bb}$, contradicting
that $r+1$ is the least exponent with this property (by definition of reduction number).
It follows that ${\rm red}_J(I)=r$ if and only if $\xx^{(r+1)\bb}\in J^{r+1}$.
Now, since $\xx^{(r+1)\bb}$ is not a multiple of an $(r+1)$-th power of any $x_i^{a_i}$ (since otherwise $\xx^{\bb}$
itself would be a multiple of that $x_i^{a_i}$), it must be the case that this membership requires the
existence of $t\geq 2$ such pure powers $x_{i_1}^{a_{i_1}},\ldots ,x_{i_t}^{a_{i_t}}$ and corresponding positive integers
$s_{i_1},\ldots,s_{i_t}$ satisfying
$$\xx^{(r+1)\bb}\in (x_{i_1}^{s_{i_1}a_{i_1}}\cdots x_{i_t}^{s_{i_t}a_{i_t}}),$$
from which our required statement follows.
\qed

\begin{Remark}\rm (1) The problem of finding the exact value of ${\rm red}_J(I)$ seems like an optimization problem.

(2) It would be interesting to make explicit $t$ above as a function of the numerical data $\{n, a_i, b_i\}$.
A tecnhicality is to clarify  in how many ways is $\xx^{(d+1)\bb}$ a multiple of a monomial generator of
the power $J^{d+1}$, where ${\rm red}_J(I)=d$.
\end{Remark}

\subsection{Binary monomials of the same degree}

An interesting driving element is to search how close the Rees ideal  is to being generated by Sylvester forms.
One might expect that if the entries are very special then the chances increase.

In this part we focus on an extreme case, where $I=(x^d,y^d,x^by^{d-b})\subset R=k[x,y]$, $d\geq 2$.
One would expect that in the case of monomials, the conditions would be wholly expressed numerically,
in terms of the given exponents.

By the general discussion in Lemma~\ref{reparametrization}, we will assume at the outset that $\gcd(d,b)=1$.

\subsubsection{The Rees algebra}

Let $\mathcal I\subset B=R[t,u,v]$ stand for the Rees ideal of
$I$ -- i.e., the kernel $\mathcal I$ of the map $t\rightarrow x^dT, u\rightarrow y^dT, v\rightarrow x^by^{d-b}T$, where $T$ is a new variable.

Let $ L\subset \mathcal I$ denote the set of generators coming from the syzygies of $I$.

We have  the following theorem.

\begin{Theorem}\label{main_binary}
With the above notation,  The Rees ideal  $\mathcal I$ is generated by $L$ and binomial Sylvester forms.
\end{Theorem}
\demo
 The starting point is the two forms in $L$, which are easily seen to be $F:=x^{d-b}v-y^{d-b}t$ and $G:=y^bv-x^bu$.
Set $d=sb+c$, with $0< c<b$.
We start the Sylvester procedure from these two, namely,  write their content $2\times 2$ matrix with respect to the ``pivot'' $(x^b,y^b)$.
 \[\left[\begin{array}{c} F \\ G \end{array}\right]=
                     \left(
    \begin{array}{cc}
      u & -v \\
      x^{d-2b}v & -y^{d-2b}t\\
    \end{array}
  \right)
                      \left[\begin{array}{c} x^b
                      \\ y^b \end{array} \right],
\]
and take its determinant $F_2:=x^{d-2b}v^2-y^{d-2b}tu$.
As observed in earlier sections, the determinant belongs to the Rees ideal $\mathcal I$.
To proceed, and provided $d-2b>b$, we take the content matrix of $F$ and $F_2$ with respect to
the same pivot $(x^b,y^b)$
and let $F_3:= x^{d-3b}v^3-y^{d-3b}tu^2$ be its determinant, which is again a Rees equation.
Continue in this fashion until $F_s:=x^{d-sb}v^s- y^{d-sb}=x^cv^s-y^ctu^{s-1}$.

Now start a new cycle, changing pivot to $(x^c, y^c)$ (since $c<b$) and write the content matrix
of $F_s$ and $F$.
If it produces a new form of $x,y$-degree $e<b$, then we take the content matrix of
$F_s$ and this latter form with pivot $(x^{\min\{c,e\}}, y^{\min\{c,e\}})$.
Proceeding in this way, along those cycles, one eventually reaches two forms
of the shape $x\TT^{\alpha_1}-y^\TT{\beta_1},\, x\TT^{\alpha_2}-y\TT^{\beta_2}$.
Taking their content with pivot $(x,y)$ yields the binomial
$\TT^{\alpha_1}\TT^{\beta_2}-\TT^{\alpha_2}\TT^{\beta_1}$,
which is the equation of the implicit toric hypersurface.

In order to make rigorous the argument along these cycles, we draw on the details of the Euclidean algorithm run on $d$ and $b$.
This idea has been employed in \cite{Dandrea} to compute generators of $\mathcal I$, but our approach will produce some simplification in order to achieve the main objective which is that of showing that the generators are Sylvester forms obtained in a sequentially ordered fashion.

Write
\begin{eqnarray}\label{Euclid}
d&=&c_1b+d_1,d_1<b\nonumber\\
b&=&c_2d_1+d_2,d_2<d_1\nonumber\\
&\vdots&\nonumber\\
d_k&=&c_{k+2}d_{k+1}+d_{k+2},d_{k+2}<d_{k+1}\\ 
&\vdots&\nonumber\\
d_{s-1}&=&c_{s+1}d_s+1, 1<d_s\nonumber\\
d_s&=&c_{s+2}\cdot 1, c_{s+2}=d_s.\nonumber
\end{eqnarray} Also denote $d_0=b$ and $d_{-1}=d$.

Introduce the following Fibonacci like sequence with moving coefficient:
$$\left\{\begin{array}{l}
e_{-1}=0, e_0=1\\
e_k=c_ke_{k-1}+e_{k-2},\, 1\leq k\leq s+2.
\end{array}
\right.$$

Fix $1\leq k\leq s+2$ and $1\leq i\leq c_k$, where $c_k$ is the $k${\rm th} successive quotient in the Euclidean algorithm.

\smallskip

We introduce the following sequences of binomials based on the above Fibonacci like sequence:
\begin{equation}\label{theFs}
F_{k,i}:=x^{d_{k-2}-id_{k-1}}v^{ie_{k-1}+e_{k-2}}-y^{d_{k-2}-id_{k-1}}t^{p_{k,i}}u^{ie_{k-1}+e_{k-2}-p_{k,i}}
\end{equation}
for $k$ odd, and
\begin{equation}\label{theGs}
G_{k,i}:=y^{d_{k-2}-id_{k-1}}v^{ie_{k-1}+e_{k-2}}-x^{d_{k-2}-id_{k-1}}t^{ie_{k-1}+e_{k-2}-q_{k,i}}u^{q_{k,i}},
\end{equation}
for $k$ even.

Here $p_{k,i}:=O_{k,i}/d$ and  $q_{k,i}:=E_{k,i}/d$), where
$O_{k,i}:=d_{k-2}-id_{k-1}+b(ie_{k-1}+e_{k-2})$ if $k$ is odd, and
$E_{k,i}:=d_{k-2}-id_{k-1}+(d-b)(ie_{k-1}+e_{k-2})$ if $k$ is even.

Note that  $F_{1,1}=F=x^{d-b}v-y^{d-b}t$, while we set $G_{0,0}:G=y^bv-x^bu$.

\smallskip

{\sc Claim 1:} $p_{k,i}$ and $q_{k,i}$ are integers satisfying
$p_{k,i}\leq ie_{k-1}+e_{k-2}$ and $q_{k,i}\leq ie_{k-1}+e_{k-2}$.

For the proof we induct on $k\geq 1$.

$k=1$:  $O_{1,i}=d_{-1}-id_0+b(ie_0+e_{-1})=d-ib+ib=d=d\cdot 1$. Of course $0\leq p_{1,i}=1\leq i=ie_0+e_{-1}$.

$k=2$:  $E_{2,i}=d_0-id_1+(d-b)(ie_1+e_0)=b-id_1+d(ic_1+1)-bic_1-b
=d(ic_1+1-i)$. Then $0\leq q_{2,i}=ic_1+1-i\leq ic_1+1=ie_1+e_0$.

Assume that $k$ is odd -- the argument for $k$ even is entirely similar. From the Euclidean Algorithm, one has $$d_{k-1}=d_{k-3}-c_{k-1}d_{k-2} \;\; {\rm and}\;\; e_{k-1}=c_{k-1}e_{k-2}+e_{k-3}.$$
Plugging these in $O_{k,i}$ one has: $$O_{k,i}= d_{k-2}(1+ic_{k-1})-id_{k-3}+b[(1+ic_{k-1})e_{k-2}+ie_{k-3}].$$
We also have by induction $$E_{k-1,1}:=d_{k-3}-d_{k-2}+(d-b)(e_{k-2}+e_{k-3})=rq_{k-1,1},$$ where $0\leq q_{k-1,1}\leq e_{k-2}+e_{k-3}$.

It obtains $$O_{k,i}=(d_{k-2}+be_{k-2})(1+ic_{k-1}-i)+ri(e_{k-2}+e_{k-3}-q_{k-1,1}).$$
Now, to conclude it remains to argue that

\smallskip

$\bullet$ If $k$ is odd, then $d$ divides $d_{k-2}+be_{k-2}$.

$\bullet$  If $k$ is even, then $d$ divides $d_{k-2}+(d-b)e_{k-2}$.

This is easily checked by induction on $k$, where $1\leq k\leq s+2$, the basic inductive relation being
$$d_{k-2}+be_{k-2}=d(\alpha-c_{k-2}\beta+c_{k-2}e_{k-3})+b\underbrace{(e_{k-2}-e_{k-4}-c_{k-2}e_{k-3})}_{0},$$
for $k$ odd, where, say, $d_{k-3}+(d-b)e_{k-3}=d\beta$ and $d_{k-4}+be_{k-4}=d\alpha$, and noting that $d_{k-4}=c_{k-2}d_{k-3}+d_{k-2}$.

(A similar inductive relation holds if $k$ is even, replacing $b$ by $d-b$)

This proves the claim.

\medskip

{\sc Claim 2:} The sequences of binomials (\ref{theFs}) and (\ref{theGs})  are Sylvester forms in an iterative way, starting out from the syzygy forms $F_{1,1}, G_{0,0}$.
In particular, they belong to $\mathcal I$.

\smallskip

We induct on $k$, separating odd and even cases.

We can start by forming the first Sylvester form stated at the beginning of the proof, using the ``pivot'' $(x^b,y^b)$ and recalling that $d_0=b$.

Let $k\geq 1$ be odd. For $i=1,\ldots,c_{k+1}-1$, suppose that $F_{k,c_k}$ and $G_{k+1,i}$ are Sylvester forms.
Write
$$\left[\begin{array}{c} F_{k,c_k}\\ G_{k+1,i}\end{array}\right]=\left(\begin{array}{cc} v^{e_k}&-t^{p_{k,c_k}} u^{e_k-p_{k,c_k}}\\ -x^{d_{k-1}-(i+1)d_k}t^{ie_k+e_{k-1}-q_{k+1,i}}u^{q_{k+1,i}}& y^{d_{k-1}-(i+1)d_k}v^{ie_k+e_{k-1}}\end{array}\right)\left[\begin{array}{c} x^{d_k}\\ y^{d_k}\end{array}\right].$$
Then, one hand,  the determinant of the $2\times 2$ matrix belongs to $\mathcal I$ and, on the other hand it equals $-G_{k+1,i+1}$ because
as $p_{k,c_k}=(d_k+be_k)/d$, one has $ie_k+e_{k-1}-q_{k+1,i}+p_{k,c_k}=(i+1)e_k+e_{k-1}-q_{k+1,i+1}$.
This proves that  $-G_{k+1,i+1}$ is a Sylvester form

Similarly, assume that $G_{k+1,c_{k+1}}$ is a Sylvester form.
Write
$$\left[\begin{array}{c} F_{k,c_k}\\ G_{k+1,c_{k+1}}\end{array}\right]=\left(\begin{array}{cc} x^{d_k-d_{k+1}}v^{e_k}&-y^{d_k-d_{k+1}}t^{p_{k,c_k}} u^{e_k-p_{k,c_k}}\\ -t^{e_{k+1}-q_{k+1,c_{k+1}}}u^{q_{k+1,c_{k+1}}}& v^{e_{k+1}}\end{array}\right)\left[\begin{array}{c} x^{d_{k+1}}\\ y^{d_{k+1}}\end{array}\right].$$
As above, $q_{k+1,c_{k+1}}=[d_{k+1}+(d-b)e_{k+1}]/d$, one obtains that the determinant of the $2\times 2$ matrix is equal to $F_{k+2,1}$.

\medskip

\medskip

Assume that $k\geq 0$ even. Supposing $c_0=0$, and keeping in mind that $d_0=b$, similarly to the odd case for $i=1,\ldots,c_{k+1}-1$ one has $$\left[\begin{array}{c} G_{k,c_k}\\ F_{k+1,i}\end{array}\right]=\left(\begin{array}{cc} -t^{e_k-q_{k,c_k}}u^{q_{k,c_k}}&v^{e_k}\\ x^{d_{k-1}-(i+1)d_k}v^{ie_k+e_{k-1}}& -y^{d_{k-1}-(i+1)d_k}t^{p_{k+1,i}}u^{ie_k+e_{k-1}-p_{k+1,i}}\end{array}\right)\left[\begin{array}{c} x^{d_k}\\ y^{d_k}\end{array}\right],$$ where the determinant of the $2\times 2$ matrix is equal to $-F_{k+1,i+1}$.

Also, $$\left[\begin{array}{c} G_{k,c_k}\\ F_{k+1,c_{k+1}}\end{array}\right]=\left(\begin{array}{cc} -x^{d_k-d_{k-1}}t^{e_k-q_{k,c_k}}u^{q_{k,c_k}}&y^{d_k-d_{k+1}}v^{e_k}\\ v^{e_{k+1}}& -t^{p_{k+1,c_{k+1}}}u^{e_{k+1}-p_{k+1,c_{k+1}}}\end{array}\right)\left[\begin{array}{c} x^{d_{k+1}}\\ y^{d_{k+1}}\end{array}\right],$$ where the determinant of the $2\times 2$ matrix is equal to $-G_{k+2,1}$.

This proves the statement of Claim 2.

\medskip

Next let
$$\Sigma:=\{G_{0,0}\}\cup\bigcup\{F_{k,i},G_{k',j}|1\leq k,k'\leq s+2,1\leq i\leq c_k,1\leq j\leq c_{k'}\}$$
stand for the set of syzygy forms and the binomial Sylvester forms obtained above.

\smallskip

{\sc Claim 3:} $\Sigma$ generates $\mathcal I$; in particular, $\mathcal I$ is generated by $1+c_1+\cdots + c_{s+2}$ syzygy forms and Sylvester forms.

\smallskip

The proof consists in showing that any binomial generator of $\mathcal I$ belongs to the ideal $(\Sigma)$.
Now, such a binomial is one of the following two types: $x^{\alpha}v^{\beta}-y^{\alpha}t^{\gamma}u^{\beta-\gamma}$ or $y^{\alpha}v^{\beta}-x^{\alpha}t^{\gamma}u^{\beta-\gamma}$.

By an apparent symmetry, it suffices to deal with the first type.
Set $H:=x^{\alpha}v^{\beta}-y^{\alpha}t^{\gamma}u^{\beta-\gamma}\in \mathcal I$.

We induct on $\alpha$.

 If $\alpha=0$, then, because $\gcd(d,b)=1$, one has $\beta=d\delta$ and $\gamma=b\delta$.
So $H=(v^d)^{\delta}-(t^bu^{d-b})^{\delta}$, hence it is divisible by the implicit equation $v^d-t^bu^{d-b}\in\Sigma$.

Thus, let $\alpha>0$ and assume the result for binomials of the two types whose $x$-degree $x$ is less than $\alpha$.

\smallskip

{\sc Claim:} There exist indices $k,i$, with $1\leq k\leq s+2$ and $1\leq i\leq c_k$, such that
$$\alpha\geq d_{k-2}-id_{k-1}\geq 1\mbox{ and }\beta\geq ie_{k-1}+e_{k-2}\geq 1.$$

Assuming the claim, one can ``divide'' $H$ by $F_{k,i}$:
$$H=x^{\alpha-(d_{k-2}-id_{k-1})}v^{\beta-(ie_{k-1}+e_{k-2})}F_{k,i}+F',$$
where $H'=x^{\alpha-(d_{k-2}-id_{k-1})}v^{\beta-(ie_{k-1}+e_{k-2})}y^{d_{k-2}-id_{k-1}}t^{p_{k,i}}u^{ie_{k-1}+e_{k-2}-p_{k,i}}- y^{\alpha}t^{\gamma}u^{\beta-\gamma}$.
Clearly, $H'$ is a multiple of a monomial in $y,t$ by a binomial $H''$ of the first type above, with $x$-degree $\alpha-(d_{k-2}-id_{k-1})<\alpha$.
Since $\mathcal I$ is a prime ideal containing no monomials, one has $H''\in \mathcal I$.
By the inductive hypothesis, $H''\in (\Sigma)$ and hence,
 $H\in (\Sigma)$.

\smallskip

In order to prove the claim we draw on the basic relation
 $$\alpha+b\beta=d\gamma$$
stemming from the assumption that $F\in \mathcal I$.

Without loss of generality we may assume that $\alpha,\beta,\gamma\geq 1$.

\medskip

\noindent\underline{First step.} Denote $\beta_1=\beta\geq 1$ and $\gamma_1=\gamma\geq 1$. We have $$\alpha+b\beta_1=d\gamma_1.$$
\begin{itemize}
  \item If $\beta_1\leq c_1$, choose $k=1$ and $i=\beta_1$. Indeed, as $\gamma_1\geq 1$, then $\alpha\geq d-\beta_1 b$
  and therefore $\alpha\geq d_{-1}-id_0\mbox{ and }\beta\geq \beta_1 e_0+e_{-1}=\beta$.

  Recall that $d_{-1}=d,d_0=b, e_{-1}=0, e_0=1$.
  \item If $\beta_1>c_1$ and $\alpha\geq d_1$, then, as $d_1=d_{-1}-c_1d_0$, choose $k=1$ and $i=c_1$.
  \item If $\beta_1>c_1$ and $\alpha<d_1$, move to the next step.
\end{itemize}

\medskip

\noindent\underline{Second step.} (Under the assumptions of the last item in the previous step).

Denote $\beta_3=\beta_1-c_1\gamma_1$ and $\gamma_3=\gamma_1-c_2\beta_3$. We still have $$\alpha+d_2\beta_3=d_1\gamma_3,$$
but need to show that $\beta_3,\gamma_3\geq 1$.

Since $d_2<d_1$ then $\alpha>(\gamma_3-\beta_3)d_1$.
If $\gamma_3-\beta_3\geq 1$, then $\alpha>d_1$ which is a contradiction. So $\gamma_3\leq \beta_3$. From $\gamma_3=\gamma_1-c_2\beta_3$, if $\beta_3<0$,
then $\gamma_3>0$ -- again, a contradiction. Also $\beta_3\neq 0$, otherwise $\alpha=d_1\gamma_3$ which contradicts $0<\alpha<d_1$. So $\beta_3>0$,
and hence from $\alpha+d_2\beta_3=d_1\gamma_3$, one gets also that $\gamma_3>0$.

\begin{itemize}
  \item If $\beta_3\leq c_3$, choose $k=3$ and $i=\beta_3$. Indeed, as $\gamma_3\geq 1$, then $\alpha\geq d_1-id_2$.
  Also $\beta\geq ie_2+e_1=\beta-c_1(\gamma_3-1),$ since $\gamma_3\geq 1$.
  \item If $\beta_3>c_3$ and $\alpha\geq d_3$, then, as $d_3=d_1-c_3d_2$, can choose $k=3$ and $i=c_3$.
  Indeed, $\beta\geq c_3e_2+e_1=c_1(c_2c_3+1)+c_3$, since $\beta=\beta_3+c_1\gamma_1>c_3+c_1\gamma_1$ and $\gamma_1=\gamma_3+c_2\beta_3>1+c_2c_3$.
  \item If $\beta_3>c_3$ and $\alpha<d_3$, move to the next step.
\end{itemize}

\medskip

\noindent\underline{General step.} (Under the assumptions of  the previous step, namely, for all $1\leq j\leq m-1$,
$\beta_{2j+1}=\beta_{2j-1}-c_{2j-1}\gamma_{2j-1}\geq 1$, $\gamma_{2j+1}=\gamma_{2j-1}-c_{2j}\beta_{2j+1}\geq 1$,
$\beta_{2j+1}>c_{2j+1}$ and $\alpha<d_{2m-1}$.)

\smallskip

Denote $\beta_{2m+1}=\beta_{2m-1}-c_{2m-1}\gamma_{2m-1}$ and $\gamma_{2m+1}=\gamma_{2m-1}-c_{2m}\beta_{2m+1}$.

\smallskip

We contend that, for every $0\leq j\leq m$, one has
\begin{equation}\label{beta}
\beta=e_{2j}\beta_{2j+1}+e_{2j-1}\gamma_{2j+1}.
\end{equation}
To see this we induct on $j$.
If $j=0$, then $\beta=e_0\beta_1+e_{-1}$. For $j>0$, have

\begin{eqnarray}
\beta &=& e_{2j-2}\beta_{2j-1}+e_{2j-3}\gamma_{2j-1}\nonumber\\
 &=&e_{2j-2}(\beta_{2j+1}+c_{2j-1}\gamma_{2j-1})+e_{2j-3}\gamma_{2j-1}\nonumber\\
 &=&e_{2j-2}\beta_{2j+1}+(c_{2j-1}e_{2j-2}+e_{2j-3})\gamma_{2j-1}\nonumber\\
 &=&e_{2j-2}\beta_{2j+1}+e_{2j-1}\gamma_{2j-1}\nonumber\\
 &=&e_{2j-2}\beta_{2j+1}+e_{2j-1}(\gamma_{2j+1}+c_{2j}\beta_{2j+1})\nonumber\\
 &=&(e_{2j-2}+c_{2j}e_{2j-1})\beta_{2j+1}+e_{2j-1}\gamma_{2j+1}\nonumber\\
 &=&e_{2j}\beta_{2j+1}+e_{2j-1}\gamma_{2j+1},\nonumber
\end{eqnarray}
 where the first equality comes from the inductive hypothesis.

We also have $$\alpha+d_{2m}\beta_{2m+1}=d_{2m-1}\gamma_{2m+1}.$$

By a similar argument to the one in the second step above (with the indices $1,2,3$ replaced by $2m-1,2m,2m+1$) we get $\beta_{2m+1},\gamma_{2m+1}\geq 1$.

\begin{itemize}
  \item If $\beta_{2m+1}\leq c_{2m+1}$, then $d_{2m-1}\gamma_{2m+1}-\alpha\leq d_{2m}c_{2m+1}=d_{2m-1}-d_{2m+1}$,
  giving $\alpha-d_{2m+1}\geq d_{2m-1}(\gamma_{2m+1}-1)$. As $\alpha<d_{2m-1}$, one gets $\gamma_{2m+1}=1$.
  Choose $k=2m+1$ and $i=\beta_{2m+1}$. Indeed, then $\alpha\geq d_{2m-1}-id_{2m}=\alpha\geq 1$.
  Also, by (\ref{beta}) one gets $\beta\geq ie_{2m}+e_{2m-1}=\beta\geq 1$.
  \item If $\beta_{2m+1}>c_{2m+1}$ and $\alpha\geq d_{2m+1}$, we take $k=2m+1$ and $i=c_{2m+1}$.
  Indeed, $\alpha\geq d_{2m-1}-id_{2m}=d_{2m+1}\geq 1$ and $\beta>c_{2m+1}e_{2m}+e_{2m-1}\geq 1$ from CLAIM above and since $\beta_{2m+1}>c_{2m+1}$.
  \item If $\beta_{2m+1}>c_{2m+1}$ and $\alpha< d_{2m+1}$ we move to step $m+1$.
\end{itemize}

Since $d_{s+1}=1, d_{s+2}=0$ and $\alpha\geq 1$, the procedure has to eventually come to a halt.
\qed

\begin{Example}\rm
Let $I=(x^{14}, y^{14},x^3y^{ 11})$. Following the prescription given by the theory so far, we find, first of all that $\mathcal I$ is generated by $1+4+1+2=8$ binomials. These may be found to be
\begin{eqnarray*} y^3v-x^3u,&&  \kern-10pt y^{11}t-x^{11}v, \;y^8tu-x^8v^2, \;y^5tu^2-x^5v^3,\; y^2tu^3-x^2v^4\,|\\
&& xtu^4-yv^5, \;yt^2u^7-xv^9,\; t^3u^{11}-v^{14}.
\end{eqnarray*}
The top line corresponds to the first cycle: $14=4\cdot 3+2$.
Note how the exponents of $x$ and $y$ in these monomials decrease steadily by ratio of $3$.
From the viewpoint of Sylvester forms, as described in (a) above, the first binomial on the bottom line is the result of pairing the last binomial on the top line and the
syzygetic binomial $y^3v-x^3u$ with pivot $(x^2,y^2)$; the one after it on the bottom line is the result of
pairing it and the one before it (of $\TT$-degree $4$)  with pivot $(x,y)$.
\end{Example}

\medskip

\subsubsection{Almost Cohen--Macaulayness}

The almost Cohen--Macaulayness of the ideals of the form $I=(x^d,y^d,x^by^{d-b})\subset R=k[x,y]$, $d\geq 2$, was established in the more general result \cite[Proposition 1.9]{RoSw} through the depth of the associated graded ring of $I$.
In this part we will derive this result directly from the nature of the defining ideal of the Rees algebra of $I$.
Moreover, we complement the result by showing that the latter is a {\em strict} almost Cohen--Macaulay if and only if $d\ge 3$.

\begin{Theorem}\label{main2}
Let $R=k[x,y]$ and $I=(x^d, y^d, x^by^{d-b})$, where $d\geq 2$ and $b\leq d-b$. Suppose that $\gcd(d,b)=1$.
Then the Rees algebra $\mathcal{R}_R(I)$ is almost Cohen--Macaulay and it is Cohen--Macaulay if and only if $d=2$.
\end{Theorem}
\demo  We induct on $0\leq k\leq s+1$, moving down from one Sylvester cycle to the next one, along the sequence $$b>d_1>d_2>\cdots>d_s>d_{s+1}=1.$$
The proof is of a telescopic nature, pretty much as the one of Proposition~\ref{resolution}: we show that $\mathcal{J}$ is obtained as a ``limit'' of subideals, each successively obtained
as a mapping cone over the previous one, and at each step the subideal has homological dimension $\leq 2$.

We introduce the telescopic subideals:
$$\mathcal T(k)=\left\lbrace
\begin{array}{cc}
(\mathcal T(k-1), F_{k,1},\ldots,F_{k,c_k}) & \mbox{if $k$ is odd}\\[6pt]
(\mathcal T(k-1), G_{k,1},\ldots,G_{k,c_k}) & \mbox{if $k$ is even}
\end{array}
\right.
$$
where  $\mathcal T(0):=(G_{0,0})$.

We may start the induction either from $k=0$, in which case the assertion is trivial, or even from $k=1$, in which case the ideal is the defining ideal $(G_{0,0}, F_{1.1})$ of the symmetric algebra of $I$, which is obviously a complete intersection.

Thus, let $k\geq 2$ be odd (to fix ideas) and assume proved that  $\mathcal T(k)$ is almost Cohen--Macaulay (meaning that $\mathcal T(k)$ has homological dimension $\leq 2$).
We wish to prove that $$\mathcal T(k+1):=(\mathcal T(k),G_{k+1,1},\ldots,G_{k+1,c_{k+1}})$$ is almost Cohen--Macaulay.

By Theorem~\ref{main_binary} and from  Lemma~\ref{mapcone} we have $$(G_{k-1,c_{k-1}},F_{k,i}):F_{k,i+1}=(x^{d_{k-1}},y^{d_{k-1}}), i=1,\ldots,c_k-1$$ and $$(G_{k-1,c_{k-1}},F_{k,c_k}):G_{k+1,1}=(x^{d_k},y^{d_k}).$$

\noindent{\bf Claim 1.} $\mathcal T(k):G_{k+1,1}=(x^{d_k},y^{d_k})$.

{\em Proof of Claim 1.} The inclusion $\supseteq$ is obvious from above. Let $f\in \mathcal T(k):G_{k+1,1}$. Then there exist $F\in (\mathcal T(k-1), F_{k,1},\ldots,F_{k,c_k-1})$ and $\beta\in R[t,u,v]$ with
$$G_{k+1,1}f=F+\beta F_{k,c_k}.$$ Multiplying this by $x^{d_k}$ one gets $$(y^{d_{k-1}-d_k}v^{e_{k-1}}F_{k,c_k}+ t^{p_{k,c_k}}u^{e_k-p_{k,c_k}}G_{k-1,c_{k-1}})f=x^{d_k}F+\beta x^{d_k}F_{k,c_k-1}.$$ One obtains this based on the fact that $e_k+e_{k-1}-q_{k+1,1}-p_{k,c_k}=q_{k-1,c_{k-1}}$ (see who are these numbers in the proof of Theorem~\ref{main_binary}).

So $y^{d_{k-1}-d_k}v^{e_{k-1}}f-\beta x^{d_k}\in (\mathcal T(k-1), F_{k,1},\ldots,F_{k,c_k-1}):F_{k,c_k}$, and the latter colon ideal is equal to the ideal $(x^{d_{k-1}},y^{d_{k-1}})$ from the previous step (i.e., step $k-1$) in the proof. Therefore $$y^{d_{k-1}-d_k}v^{e_{k-1}}f-\beta x^{d_k}=\delta_1 x^{d_{k-1}}+\delta_2 y^{d_{k-1}},$$ for some $\delta_1,\delta_2\in R[t,u,v]$.

This can be written as $$y^{d_{k-1}-d_k}(v^{e_{k-1}}f-\delta_2y^{d_k})=x^{d_k}(\delta_1x^{d_{k-1}-d_k}-\beta),$$ which leads to $x^{d_k}$ dividing $v^{e_{k-1}}f-\delta_2y^{d_k}$. So $$f\in (x^{d_k},y^{d_k}):v^{e_{k-1}}=(x^{d_k},y^{d_k}),$$ and Claim 1 is shown.

\noindent{\bf Claim 2.} $(\mathcal T(k),G_{k+1,1},\ldots,G_{k+1,j-1}):G_{k+1,j}=(x^{d_k},y^{d_k})$, for $j=1,\ldots, c_{k+1}$.

Here $G_{k+1,0}$ is to be taken as the $0$ polynomial.

{\em Proof of Claim 2.}
By Theorem~\ref{main_binary} and
Lemma\ref{mapcone}, one has  $(F_{k,c_k},G_{k+1,j-1}):G_{k+1,j}=(x^{d_k},y^{d_k})$, showing the inclusion $\supseteq$.

The reverse inclusion will be shown by induction on $j=1,\ldots,c_{k+1}$.

For $j=1$ it is clear from Claim 1.

Suppose $j\geq 2$ and let $f\in (\mathcal T(k),G_{k+1,1},\ldots,G_{k+1,j-1}):G_{k+1,j}$.
Then, there exists $G\in (\mathcal T(k),G_{k+1,1},\ldots,G_{k+1,j-2})$ and $\beta\in R[t,u,v]$ such that $$G_{k+1,j}f=G+\beta G_{k+1,j-1}.$$ Multiplying this equation by $y^{d_k}$, after appropriate substitutions one gets $$(v^{e_k}G_{k+1,j-1}+x^{d_{k-1}-jd_k}t^{*}u^{**}F_{k,c_k})f=y^{d_k}G+\beta y^{d_k}G_{k+1,j-1},$$ where $*$ and $**$ are suitable positive integers.

It yields $v^{e_k}f-\beta y^{d_k}\in (\mathcal T(k),G_{k+1,1},\ldots,G_{k+1,j-2}): G_{k+1,j-1}$, which by induction equals $(x^{d_k},y^{d_k})$. Thus
$$f\in (x^{d_k},y^{d_k}):v^{e_k}=(x^{d_k},y^{d_k}),$$
proving the claim.

Now apply recursively the mapping cone construction to deduce that $\mathcal T(k+1)$ is almost Cohen--Macaulay.

\medskip

For the complementary assertion, by \cite[Theorem 2.1]{syl1} (see also \cite[Corollary 4.2 (i)]{syl3}) $\mathcal{R}_R(I)$ is Cohen--Macaulay iff the reduction number of $I$ equals to $1$. Since $\gcd(b,d)=1$ (hence the map is birational), the implicit equation has degree $d$. Therefore, the reduction number is $d-1$. This implies $d=2$ as asserted. \qed

\subsubsection{The Huckaba--Marley approach}

An alternative to proving almost Cohen--Macaulayness goes via the Huckaba--Marley test in terms of the lengths of suitable modules (\cite[Theorem 4.7]{HM}).
Since monomial almost complete intersection  zero-dimensional ideals are Ratliff--Rush closed by \cite[Proposition 1.9]{RoSw}, then the sought result follows from \cite[Corollary 4.13 (1)]{HM}.

Since this work is about being very explicit,  we will explain the nature of  the modules used in the Huckaba--Marley approach, as well as their lengths and the relation to their annihilators.

To make sense, we harmlessly pass to the local ring $k[x,y]_{\fm}$, where $\fm=(x,y)$, and consider the extended ideal $I_{\fm}$.
We update the notation accordingly. It goes without saying that we often argue in the original
standard graded polynomial ring $k[x,y]$, then pass to its irrelevant localization.

Assuming that $\gcd(b,d)=1$ (hence the map is birational), we have $e_1(I)=e_1((x,y)^d)$, where $e_1(\_)$ denotes the second Hilbert coefficient. Therefore, $e_1(I)={{d}\choose {2}}$ (\cite[Proposition 3.3]{syl2}).

The reduction number being $d-1$, the Huckaba--Marley test tells us that  almost Cohen--Macaulayness means the inequality
$$\sum_{\ell=1}^{d-1} \lambda(I^\ell/JI^{\ell-1})\leq {{d}\choose {2}},$$ where $J=(x^d,y^d)$ and $\lambda$ denotes length.

The following result gives a general expression for the lengths above.

\begin{Lemma}
\label{lengths_are_cis} Suppose $\gcd(d,b)=1$.
For all $\ell\geq 1$ one has
$$JI^{\ell-1}:I^{\ell}=(x^{s_\ell},y^{t_\ell}),$$ for suitable integer exponents $s_\ell, t_\ell$.
Moreover, these exponents are effectively computable {\rm (}optimization{\rm )}.
\end{Lemma}
\demo
It is apparent that, for any $\ell\geq 1$, one has $I^\ell=(JI^{\ell-1}, x^{b\ell}y^{(d-b)\ell})$.
Therefore, $ \lambda(I^\ell/JI^{\ell-1})=\lambda(R/(JI^{\ell-1}:x^{b\ell}y^{(d-b)\ell}))$.

We proceed to compute the latter.

First, consider all $3$-partitions of $\ell-1$ in nonnegative integers: $\ell-1=i+j+k$.
A straightforward inspection gives that the typical generator of $JI^{\ell-1}$ has the form
$$x^{id+jb}y^{k(d+1)+j(d-b)}\mbox{ or }x^{(i+1)d+jb}y^{kd+j(d-b)}.$$
Recall the following general fact: if $m_1,\ldots,m_p$ are monomials, then $(m_1,\ldots,m_{p-1}):m_p$ is a monomial ideal generated by the positive powers of the variables in $m_i/m_p,i=1,\ldots,p-1$.
Therefore, upon replacing $k=\ell-1-(i+j)$, the generators of $JI^{\ell-1}:x^{b\ell}y^{(d-b)\ell}$ are of the form
$$x^{id-(\ell-j)b}y^{(\ell-j)b-id}\mbox{ or }x^{(i+1)d-(\ell-j)b}y^{(\ell-j)b-(i+1)d},$$ for $0\leq i+j\leq \ell-1$.
Clearly, the above monomials are of the form $x^{\varepsilon}y^{-\varepsilon}$, for some integer $\varepsilon$ and its opposite. Disregarding the negative power, this shows
that indeed $JI^{\ell-1}:x^{b\ell}y^{(d-b)\ell}=(x^{s_\ell},y^{t_\ell})$, where $s_{\ell}$ and $t_{\ell}$ are suitable exponents.

To express these exponents in further detail, we proceed as follows.
Quite clearly, observe that for  nonnegative integers $\alpha,\beta$, one has $\alpha d-\beta b>0$   if and only if $\alpha\geq 1$ and $\beta/\alpha\leq c_1$, where $d=c_1b+d_1$ is the first step in the Euclidean Algorithm.

Set $\Delta_{\ell-1}:=\{(i,j)\,\vert \, i+j\leq\ell-1,i\geq 0,j\geq 0\}$. By the latter observation, $\Delta_{\ell-1}$ can be written as the disjoint union of the following three sets
\begin{enumerate}
  \item $\Gamma_{\ell}:=\{(i,j)\in\Delta_{\ell-1}|i\geq 1,ic_1+j\geq\ell\}$;
  \item $\Omega_{\ell}:=\{(i,j)\in\Delta_{\ell-1}|i\geq 1,ic_1+j<\ell\leq(i+1)c_1+j\}\cup\{(0,j)\in\Delta_{\ell-1} |c_1+j\geq\ell\}$;
  \item $\Lambda_{\ell}:=\{(i,j)\in\Delta_{\ell-1}|(i+1)c_1+j<\ell\}$,
\end{enumerate}
based on the conditions
\begin{enumerate}
  \item $id-(\ell-j)b>0$;
  \item $id-(\ell-j)b<0$ and $(i+1)d-(\ell-j)b>0$;
  \item $(i+1)d-(\ell-j)b<0$.
\end{enumerate}

Consider the following numbers

\begin{eqnarray}\displaystyle
m_{\ell}&=&\min_{(i,j)\in \Gamma_{\ell}}\{id-(\ell-j)b\}\nonumber\\
m_{\ell}'&=&\min_{(i,j)\in \Omega_{\ell}}\{(i+1)d-(\ell-j)b\}\nonumber\\
n_{\ell}'&=&\min_{(i,j)\in \Omega_{\ell}}\{(\ell-j)b-id\}\nonumber\\
n_{\ell}&=&\min_{(i,j)\in \Lambda_{\ell}}\{(\ell-j)b-(i+1)d\}.\nonumber
\end{eqnarray}

Then
$$s_{\ell} =
\left\{
\begin{array}{ll}
 m_{\ell}, & \mbox{\rm if $\Gamma_{\ell}\neq \emptyset$}\\
  m'_\ell, & \mbox{\rm  otherwise.}\\
\end{array}
\right.
$$
and, similarly
$$t_{\ell} =
\left\{
\begin{array}{ll}
 n_{\ell}, & \mbox{\rm  if $\Lambda_{\ell}\neq \emptyset$}\\
  n'_\ell, & \mbox{\rm  otherwise.}\\
\end{array}
\right.
$$

\qed

\begin{Example}\rm Suppose we are in the particular case where that $d=cb+1, b\geq 2$. Let us assume further that $c\geq 2$. Then we have the following:
\end{Example}
\begin{itemize}
  \item $\ell=1$. Then $\Gamma_1=\Lambda_1=\emptyset$, hence $m_1$ and $n_1$ do not exist, and $m_1'=d-b,n_1'=b$, giving $s_1=d-b,t_1=b$.
  \item $2\leq\ell\leq c$. $\Lambda_\ell=\emptyset$, hence $n_{\ell}$ does not exist. Simple calculations will show that $m_{\ell}=m_{\ell}'=d-\ell b$, and $n_{\ell}'=b$, giving that $s_{\ell}=d-\ell b, t_{\ell}=b$.
  \item $\ell\geq c+1$. In this case $(1,\ell-c)\in\Gamma_{\ell}$, and since $1\cdot d-(\ell-(\ell-c))b=1$, one gets $m_{\ell}=1$ -- smallest possible, hence $s_{\ell}=1$.

Detecting $t_{\ell}$ is harder when $\ell\geq c+1$.
Suppose that $c+1\leq\ell\leq 2c$. Then $\Omega_{\ell}=\{(1,0),\ldots,(1,\ell-c-1)\}\cup\{(0,\ell-c),\ldots,(0,\ell-1)\}$, and $\Lambda_{\ell}=\{(0,0),\ldots,(0,\ell-c-1)\}$. In this case $n_{\ell}'=n_{\ell}=b-1$, hence $t_{\ell}=b-1$.

More generally, let
$\delta c+1\leq\ell\leq (\delta+1)c$, for $\delta\in\{2,\ldots,b-1\}$. Observe that for $\delta=b-1$ we have $(\delta+1)c=bc=d-1$ the extreme value of $\ell$.
For arbitrary $\delta$ one has $(\delta,\ell-\delta c-1)\in \Omega_{\ell}$, and hence
$$t_{\ell}\leq [\ell-(\ell-\delta c-1)]b-\delta d=\delta(cb-d)+b=b-\delta.$$
\end{itemize}
(In fact, equality holds, but a direct calculation is harder.)

In any case, we have enough information to derive the following inequality:

\begin{eqnarray}\displaystyle
\sum_{\ell=1}^{d-1}s_{\ell}t_{\ell}&=&\sum_{\ell=1}^cb(d-\ell b)+\sum_{\ell=c+1}^{d-1}1\cdot t_{\ell}\nonumber\\
&\leq&\sum_{\ell=1}^cb(d-\ell b)+\sum_{\delta=1}^{b-1}c(b-\delta)\nonumber\\
&=&cbd-\frac{c(c+1)}{2}b^2+cb(b-1)-c\frac{b(b-1)}{2}\nonumber\\
&=&\frac{cb}{2}\left[2d-(c+1)b+2(b-1)-(b-1)\right]\nonumber\\
&=&\frac{cb}{2}\left[2d-cb-1\right]\nonumber\\
&=&\frac{(d-1)d}{2}={{d}\choose{2}},\nonumber
\end{eqnarray} as $d=cb+1$.

\medskip

We next elaborate  on the previous data in the general case.

For $\ell=1,2$ the above exponents are within reach.
They can also be obtained directly, as follows.
Fo $\ell=1$, one has $J:I=(x^b,y^{d-b})$ as it follows directly from the syzygy matrix of $I$.
For $\ell=2$  one can use \cite[Proposition]{syl3}, by which $\lambda(I^2/JI)=\lambda(I/J)-\lambda (R/I_1(\phi))$. Since $b\leq d-b$ by assumption, then $\lambda(R/I_1(\phi))=b^2$ so it follows that $\lambda(I^2/JI)= b(d-2b)$.

Thus, together the first two lengths contribute $b(2d-3b)$.

Moreover, one has:

\begin{Lemma}
\label{linear syzygies}
Notation as in the previous lemma.
Then:
\begin{enumerate}
\item[{\rm (i)}] There is a least index $\ell_0$ for which either $s_{\ell_0}=1$ or $t_{\ell_0}=1$
\item[{\rm (ii)}] There is a least index $\ell'_0$ for which both $s_{\ell_0}=t_{\ell_0}=1$
\item[{\rm (iii)}] If $\ell_0$ and $\ell'_0$ are as in {\rm (i)} and {\rm (ii)}, then $\ell'_0\geq d-\ell_0$.
\end{enumerate}
\end{Lemma}
\demo We prove (ii) -- note that (i) is a  consequence thereof.
The equality $JI^{\ell-1}:I^{\ell}=(x,y)$ means that $I^{\ell}$ has two independent linear syzygies, which is tantaumont to having two forms in $\mathcal I$ of bidegree $(1,\ell)$.
Saying that $\ell$ is the first index for such an occurrence implies that they are independent in this bidegree.
The existence of such equations is forced by birationality (which we are assuming with $\gcd(d,b)=1$) according to the main criterion in  \cite[Theorem 2.18]{aha}.

(iii) This follows again from the further details of \cite[Theorem 2.18]{aha}: the Jacobian matrix of all such equations as in (ii) with respect to $x,y$ has to have rank $1$ modulo the implicit equation.
Since the latter has degree $d$, the $2$-minors of this matrix must be multiples of the implicit equation, hence of degree at least $d$.
\qed

\begin{Question}\rm
{\bf (1)} Are the distinguished indices  $\ell_0$ and $\ell_0'$ equidistant from the extremes
of the sequence $\{1, 2,\ldots, d-2,d-1\}$, i.e., is $\ell_0'=d-\ell_0$ always?

{\bf (2)}  The length sequence is divided into three basic sectors:
(i) from $\ell=1$ through $\ell_0-1$; (ii) from $\ell_0$ through $\ell_0'-1$; (iii) from $\ell_0'$ on.
The sector (iii) contributes length $1+\cdots +1=d-1-\ell_0'+1=d-\ell_0'$; the sector (ii) comes from quotients (say)
$JI^{\ell-1}:I^{\ell}= (x,y^{t_{\ell}})$,
thus contributes length $\sum _{\ell_0}^{\ell_0'-1} t_{\ell}$ (is there a pattern for $t_{\ell}$?).
The sector (i) has the big summand $b(2d-3b)$, but otherwise does not exhibit any particular pattern for the subsequent lengths.
Also note that, in contrast to the case of general forms, here even $I^2$ can have a linear syzygy
(e.g., if $b=3, d=7$).
\end{Question}

\subsection{The uniform case}

By {\em uniform} we mean the case where $a_1=\cdots =a_n:=a$ and $b_1=\cdots =b_n:=b$, with $0<b<a$.

The driving force in this part is the following conjecture:

\begin{Conjecture}\label{conjecture_uniform}
Let $I=(x_1^a,\ldots, x_n^a, (x_1\cdots x_n)^b),$ where $0<b<a$.
Then $\mathcal{R}_R(I)$ is almost Cohen--Macaulay and is Cohen--Macaulay if and only if $a\leq 2b$.
\end{Conjecture}

The clear part of the statement so far is the implication $a\leq 2b \Rightarrow \mathcal{R}_R(I)$ is Cohen--Macaulay.
Indeed, the inequality implies a quadratic equation arising from the inclusion $(x_1\cdots x_n)^{2b}\in (x_1^a\cdots x_n^a)$, hence
$I^2\subset JI$, where $J=(x_1^a,\ldots, x_n^a)$.
Therefore $J$ is a minimal reduction with reduction number $1$.

\medskip

We next state some further details about reduction numbers in this situation.

\begin{Proposition}\label{red_no_ab} Let $I=(x_1^a,\ldots,x_n^a,(x_1\cdots x_n)^b)$, with $0<b<a$. The following hold:
\begin{itemize}
\item[{\rm (a)}] $J:=(x_1^a,\ldots,x_n^a)$ is a minimal reduction of $I$ if and only if $nb\geq a${\rm ;} in this case, letting $1\leq p\leq n$ be the smallest integer such that $pb\geq a$ {\rm (}hence $(p-1)b<a${\rm )}, one has ${\rm red}_J(I)=p-1$.
\item[{\rm (b)}] If $nb<a$, then
$Q:=(x_1^a-x_n^a,\ldots,x_{n-1}^a-x_n^a,(x_1\cdots x_n)^b)$ is a minimal reduction of $I$ and ${\rm red}_Q(I)=n-1$.
\end{itemize}
\end{Proposition}
\demo (a)
Suppose that $J$ is a minimal reduction, and let ${\rm red}_J(I)=r$. Then, by Lemma \ref{red_no}, there exist $n\geq t\geq 2$
and $s_{i_1},\ldots,s_{i_t}$ with $s_{i_1}+\cdots+s_{i_t}=r+1$ such that
$$(r+1)b\geq s_{i_j}a, j=1,\ldots,t.$$
Adding up the inequalities one gets $tb\geq a$ and hence, $nb\geq a$.

Conversely, letting $J:=(x_1^a,\ldots,x_n^a)$, since $$((x_1\cdots x_n)^b)^p\in(x_1^a\cdots x_n^a),$$ one obtains that $JI^{p-1}=I^p$, and hence ${\rm red}_J(I)\leq p-1$.
Suppose that ${\rm red}_J(I)= p-q, q\geq 2$. Then, by Lemma \ref{red_no}, there exist at least one $1\leq\ell\leq t$, such that
$$(p-q+1)b\geq s_{i_{\ell}}a.$$
This is a contradiction, since $a>(p-1)b\geq(p-q+1)b$, and $s_{i_{\ell}}a\geq a$.

\medskip

(b)
Se $Q=(x_1^a-x_n^a,\ldots,x_{n-1}^a-x_n^a,(x_1\cdots x_n)^b)$.
We first claim that $I^n\subset QI^{n-1}$.

Thus, let $$\mathcal M=x_1^{i_1a}\cdots x_n^{i_na}(x_1\cdots x_n)^{bj},i_1+\cdots+i_n+j=n$$
be a typical generator of $I^n$.

Suppose that for some $1\leq s\leq n-1$, $i_s\geq 1$. Then
$$x_s^{i_sa}=x_s^{(i_s-1)a}x_s^a = \underbrace{x_s^{(i_s-1)a}(x_s^a-x_n^a)}_{\in Q}+x_s^{(i_s-1)a}x_n^a.$$
We get that $\mathcal M=\mathcal M'+\mathcal M'',$ where $\mathcal M'\in QI^{n-1}$, and
$$\mathcal M''=x_1^{i_1a}\cdots x_s^{(i_s-1)a}\cdots x_{n-1}^{i_{n-1}a} x_n^{(i_n+1)a}(x_1\cdots x_n)^{bj}.$$
Of course, $\mathcal M\in QI^{n-1}$ iff $\mathcal M''\in QI^{n-1}$.

Repeating the process we derive that $\mathcal M\in QI^{n-1}$ exactly when $\mathcal N:=x_n^{(n-j)a}(x_1\cdots x_n)^{bj}\in QI^{n-1}$.

If $j>0$, then
$$\mathcal N=\underbrace{(x_1\cdots x_n)^b}_{\in Q}\underbrace{(x_n^a)^{(n-j)}((x_1\cdots x_n)^b)^{(j-1)}}_{\in I^{n-1}}.$$

If $j=0$, then $\mathcal N=x_n^{na}$. Using the generators $x_i^a-x_n^a\in Q, 1\leq i\leq n-1$, we have that $\mathcal N\in QI^{n-1}$
if and only if $x_1^a\cdots x_n^a\in QI^{n-1}$. But the latter is always the case because
$$x_1^a\cdots x_n^a=\underbrace{(x_1\cdots x_n)^b}_{\in Q}\underbrace{(x_1\cdots x_n)^{a-b}}_{\in I^{n-1}},$$
as $a-b>(n-1)b$.

\medskip

To complete the proof, we have to show that $I^{n-1}\not\subset QI^{n-2}$.  Since $x_n^{(n-1)a}\in I^{n-1}$, it is enough to show that $x_n^{(n-1)a}\notin QI^{n-2}$.
Suppose the contrary. Then
$$x_n^{(n-1)a}= \sum C_{(i_1,\ldots,i_n)}^{k,j}(x_k^a-x_n^a)x_1^{i_1a+bj}\cdots x_n^{i_na+bj}+ \sum P_{(i_1,\ldots,i_n)}^jx_1^{i_1a+b(j+1)}\cdots x_n^{i_na+b(j+1)},$$
where the sums are taken over all $1\leq k\leq n-1$ and $i_1+\cdots+i_n+j=n-2$.

First, observe that when $j=0$, $C_{(i_1,\ldots,i_n)}^{k,0}$ are constant polynomials.

The terms which are pure powers of $x_n$ in the righthand side
have $j=i_1=\cdots=i_{n-1}=0, i_n=n-2$. It follows that
$$x_n^{(n-1)a}=(-\sum_k C_{(0,\ldots,0,n-2)}^{k,0})x_n^{(n-1)a}+ \sum_k C_{(0,\ldots,0,n-2)}^{k,0}x_k^ax_n^{(n-2)a}+\cdots.$$
Hence $-\sum_k C_{(0,\ldots,0,n-2)}^{k,0}=1$ and the coefficients of all the other monomials must be zero.

The monomial $x_k^ax_n^{(n-2)a}$ also can occur only in $(x_k^a-x_n^a)x_k^ax_n^{(n-3)a}$. Therefore we have
$$0=\sum_k(C_{(0,\ldots,0,n-2)}^{k,0}-C_{(0,\ldots,1,\ldots,0,n-3)}^{k,0})x_k^ax_n^{(n-2)a}+ \sum_k C_{(0,\ldots,1,\ldots,0,n-3)}^{k,0}x_k^{2a}x_n^{(n-3)a}+\cdots.$$
The $1$ in the multi-index above occurs in position $k$.

We get that $C_{(0,\ldots,0,n-2)}^{k,0}-C_{(0,\ldots,1,\ldots,0,n-3)}^{k,0}=0$ for all $1\leq k\leq n-1$. If we repeat the process in the end we obtain
$$C_{(0,\ldots,0,n-2)}^{k,0}=C_{(0,\ldots,1,\ldots,0,n-3)}^{k,0}= \cdots=C_{(0,\ldots,n-2,\ldots,0,0)}^{k,0},$$
and
$$0=\sum_kC_{(0,\ldots,n-2,\ldots,0,0)}^{k,0}x_k^{(n-2)a}+\mbox{ other terms not pure powers of the variables}.$$
This leads to $C_{(0,\ldots,0,n-2)}^{k,0} = C_{(0,\ldots,n-2,\ldots,0,0)}^{k,0}=0$ for all $k$.
But this contradicts the fact that $\sum_k C_{(0,\ldots,0,n-2)}^{k,0}=-1$.
\qed

\subsection{The case of three variables}

We observe that the uniform case in two variables is uninteresting because the reduction number is $1$, hence the Rees algebra is Cohen--Macaulay.
In this subsection we focus on the case of three variables, which is already quite significant.

 For easier reading, we change the notation: $I=(x^a,y^a,z^a,(xyz)^b)\subset R:=k[x,y,z]$.
 By the obvious part of Conjecture~\ref{conjecture_uniform} (see the commentary right after it), we assume that $a>2b$.
 By Proposition~\ref{red_no_ab}, the reduction number of $I$ is $r=3-1=2$ in any case.

 Denote $t,u,v,w$  the ``external'' variables corresponding orderly to the generators of $I$.
 Set $S:=R[t,u,v,w]$.

 Let $\mathcal I\subset S$ denote the presentation ideal of the Rees algebra of $I$ on $S$ and let $ L\subset \mathcal I$ denote the set of generators coming from the syzygies of $I$.

 The following is the main result of this part.

\begin{Theorem}\label{main3} With the above notation one has:
\begin{enumerate}
\item[{\rm (a)}] The presentation ideal of the Rees algebra $\mathcal{R}_R(I)$ is generated by $L$ and by Sylvester forms originated from it.
\item[{\rm (b)}] $\mathcal{R}_R(I)$  is almost Cohen--Macaulay.
\end{enumerate}
\end{Theorem}
\demo
(a) Assume that $a\leq 3b$.

The syzygy module of $I$ is generated by the columns of the following $4\times 6$ matrix
$$\left(
  \begin{array}{cccccc}
    -y^a & -z^a & 0 & -(yz)^b & 0 & 0\\
    x^a  &   0  & -z^a & 0 & -(xz)^b & 0\\
    0    &  x^a & y^a & 0 & 0 & -(xy)^b\\
    0    &   0  &  0  & x^{a-b} & y^{a-b} & z^{a-b}
  \end{array}
\right).
$$

From this we get the following syzygy forms in the Rees ideal $\mathcal I\subset S:=R[t,u,v,w]$
\begin{eqnarray} L=
& &f_1:=x^au-y^at;\, f_2:=x^av-z^at;\, f_3:=y^av-z^au;\nonumber\\
& &g_1:=x^{a-b}w-(yz)^bt;\, g_2:=y^{a-b}w-(xz)^bu g_3:=z^{a-b}w-(xy)^bv;\nonumber
\end{eqnarray}

We form the Sylvester content matrix of $\{g_1,g_2\}$ with respect to the complete intersection $\{x^b,y^b\}$:

$$\left[\begin{array}{c}g_1\\g_2\end{array}\right]=\underbrace{\left(\begin{array}{cc}x^{a-2b}w & -z^bt\\-z^bu& y^{a-2b}w \end{array}\right)}_{\mathcal M_{12}}\left[\begin{array}{c}x^b\\y^b\end{array}\right].$$
Set $H_1:=\det(\mathcal M_{12})=(xy)^{a-2b}w^2-z^{2b}tu$.

Similarly, for  $\{g_1,g_3\}$ and $\{g_2,g_3\}$ with respect to $\{x^b,z^b\}$ and $\{y^b,z^b\}$, respectively:

$$\left[\begin{array}{c}g_1\\g_3\end{array}\right]=\underbrace{\left(\begin{array}{cc}x^{a-2b}w & -y^bt\\-y^bv& z^{a-2b}w \end{array}\right)}_{\mathcal M_{13}}\left[\begin{array}{c}x^b\\z^b\end{array}\right].$$
Set $H_2:=\det(\mathcal M_{13})=(xz)^{a-2b}w^2-y^{2b}tv$.

$$\left[\begin{array}{c}g_2\\g_3\end{array}\right]=\underbrace{\left(\begin{array}{cc}y^{a-2b}w & -x^bu\\-x^bv& z^{a-2b}w \end{array}\right)}_{\mathcal M_{23}}\left[\begin{array}{c}y^b\\z^b\end{array}\right].$$
Set $H_3:=\det(\mathcal M_{23})=(yz)^{a-2b}w^2-x^{2b}uv$.

Next take the Sylvester form of  $\{g_1,H_3\}$ with respect to $\{x^{a-b},(yz)^{a-2b}\}$:

$$\left[\begin{array}{c}g_1\\H_3\end{array}\right]=\underbrace{\left(\begin{array}{cc}w & -(yz)^{3b-a}t\\-x^{3b-a}uv& w^2 \end{array}\right)}_{\mathcal N_{13}}\left[\begin{array}{c}x^{a-b}\\(yz)^{a-2b}\end{array}\right].$$
Set $E:=\det(\mathcal N_{13})=w^3-(xyz)^{3b-a}tuv$.

Repeating the analogous procedure for $\{g_2,H_2\}$ with respect to $\{y^{a-b},(xz)^{a-2b}\}$ and of  $\{g_3,H_1\}$ with respect to $\{z^{a-b},(xy)^{a-2b}\}$ will give us again $E$.

\smallskip

If $a>3b$, then the only change is that we get $E':=(xyz)^{a-3b}w^3-tuv$ instead, and this can be obtained in three ways, similar as above, except that we replace the vector $\left[\begin{array}{c}x^{a-b}\\(yz)^{a-2b}\end{array}\right]$ by $\left[\begin{array}{c}x^{2b}\\(yz)^{b}\end{array}\right]$ (and the corresponding permutations of the variables).

\medskip

{\sc Claim:} $\mathcal I=(L, H_1,H_2,H_3,E)$.

\smallskip

The strategy is a greedy scheme showing that any binomial $H:=M_1-M_2\in \mathcal I$ is contained in the ideal of the right hand side in the claim.
Moreover, since $\mathcal I$ is a prime ideal, one may assume that $\gcd\{M_1,M_2\}=1$. As  $\mathcal I$ is generated by such binomials, we will be through.

We can assume at the outset that $H$ effectively involves the variable $w$ as otherwise it would be a relation of $(x^a,y^a,z^a)$, and hence  $H\in (L)$.

 Assuming such a binomial $H$, let $\delta,\alpha_1,\alpha_2,\alpha_3,\beta_1,\beta_2,\beta_3\geq 0$ denote the degrees in $H$ of the variables $w,t,u,v,x,y, z$, respectively.

First we have $0\leq\beta_i<a$, otherwise if, e.g., $\beta_1\geq a$, $x^{\beta_1}$ is replaced by $x^{\beta_1-a}t$. Further, if for example $w^{\delta}t^{\alpha_1}|M_1, \alpha_1\geq 1$, then we must have $x^{\beta_1}|M_2$ and $a\leq \delta b+\alpha_1 a=\beta_1$, contradiction. This is saying that if $w|M_1$, then none of $t,u$ or $v$ can divide $M_1$, this leading also to $\alpha_1+\alpha_2+\alpha_3=\delta$.

This way, up to permutation of the variables, we may assume that $H$ belongs to one of the following types:

\begin{itemize}
  \item Type 1: $w^{\delta}-x^{\beta_1}y^{\beta_2}z^{\beta_3}t^{\alpha_1}u^{\alpha_2}v^{\alpha_3}$;
  \item Type 2: $x^{\beta_1}w^{\delta}-y^{\beta_2}z^{\beta_3}t^{\alpha_1}u^{\alpha_2}v^{\alpha_3}, \beta_1>0$;
  \item Type 3: $x^{\beta_1}y^{\beta_2}w^{\delta}-z^{\beta_3}t^{\alpha_1}u^{\alpha_2}v^{\alpha_3}, \beta_1,\beta_2>0$;
  \item Type 4: $x^{\beta_1}y^{\beta_2}z^{\beta_3}w^{\delta}-t^{\alpha_1}u^{\alpha_2}v^{\alpha_3}, \beta_1,\beta_2,\beta_3>0$.
\end{itemize}

{\sc Claim:} One can assume that $\delta\leq 3$.

Supposing that $\delta\geq 4$, we consider each one of the above listes types.

\medskip

$\bullet$ Let $H$ be of Type 1. If $\delta\geq 4$, then $\alpha_i\geq 1, i=1,2,3$, otherwise we'd obtain for example $\beta_1\geq 4b>a$. Then we can write $$H=w^{\delta-3}(w^3-(xyz)^{3b-a}tuv)+w^{\delta-3}(xyz)^{3b-a}tuv-x^{\beta_1}y^{\beta_2}z^{\beta_3}t^{\alpha_1} u^{\alpha_2}v^{\alpha_3}$$ $$=w^{\delta-3}(w^3-(xyz)^{3b-a}tuv)+tuv(w^{\delta-3}(xyz)^{3b-a}-x^{\beta_1}y^{\beta_2}z^{\beta_3}t^{\alpha_1-1} u^{\alpha_2-1}v^{\alpha_3-1}).$$
Since $\mathcal I$ is a prime ideal, $w^{\delta-3}(xyz)^{3b-a}-x^{\beta_1}y^{\beta_2}z^{\beta_3}t^{\alpha_1-1} u^{\alpha_2-1}v^{\alpha_3-1}\in \mathcal I$, and hence we induct on $\delta$.

\smallskip

$\bullet$ Let $H$ be of Type 2. If $3b\geq a$, then replace $w^{\delta}$ in $H$ with $w^{\delta-3}(w^3-(xyz)^{3b-a}tuv+(xyz)^{3b-a}tuv)$ and proceed as in the previous type.

Suppose $3b<a$. If $\beta_1\geq a-b$, then replace $x^{\beta_1}w^{\delta}$ by $x^{\beta_1-(a-b)}w^{\delta-1}(x^{a-b}w-(yz)^bt+(yz)^bt)$ and apply the previous procedure.

If $\beta_1<a-b$, then $\alpha_2$ and $\alpha_3$ cannot both be zero, otherwise $\alpha_1=\delta$ and hence $\beta_1=\delta(a-b)\nless a-b$. If $\alpha_2,\alpha_3\geq 1$, since $\alpha_1\geq 1$, then replace $t^{\alpha_1}u^{\alpha_2}v^{\alpha_3}$ in $H$ with $t^{\alpha_1-1}u^{\alpha_2-1}v^{\alpha_3-1}(tuv-(xyz)^{a-3b}w^3+(xyz)^{a-3b}w^3)$ and use the same trick again.
If $\alpha_2\geq 1$ and $\alpha_3=0$, then $\beta_3=\delta b>2b$, and replace $z^{\beta_3}t^{\alpha_1}u^{\alpha_2}v^{\alpha_3}$ with $z^{\beta_3-2b}t^{\alpha_1-1}u^{\alpha_2-1}v^{\alpha_3}(z^{2b}tu-(xy)^{a-2b}w^2+(xy)^{a-2b}w^2)$, and conclude as before.

\smallskip

$\bullet$ Let $H$ be of Type 3. If $3b\geq a$, then replace $w^{\delta}$ in $H$ with $w^{\delta-3}(w^3-(xyz)^{3b-a}tuv+(xyz)^{3b-a}tuv)$ and repeat the same act.

Suppose $3b<a$. We have that $\alpha_1,\alpha_2\geq 1$. Depending on $\alpha_3=0$ or $\alpha_3\geq 1$ we use a similar argument as in the previous passages.

$\bullet$ Let $H$ be of Type 4. Then $a\geq 3b$ and $\alpha_i\geq 1,i=1,2,3$, so then replace $t^{\alpha_1}u^{\alpha_2}v^{\alpha_3}$ in $H$ with $t^{\alpha_1-1}u^{\alpha_2-1}v^{\alpha_3-1}(tuv-(xyz)^{a-3b}w^3+(xyz)^{a-3b}w^3)$ and we implement the same downgrading procedure.

\medskip

Thus, we assume that $\delta\leq 3$.
In this case we will show that $H$ is actually one of the binomials on the right hand side ideal of the statement.

\medskip

\noindent \underline{Type 1:} $H=w^{\delta}-x^{\beta_1}y^{\beta_2}z^{\beta_3}t^{\alpha_1}u^{\alpha_2}v^{\alpha_3}$. Then we have

\begin{eqnarray}
\delta b&=&\beta_1+\alpha_1a \nonumber \\
\delta b&=&\beta_2+\alpha_2a \nonumber\\
\delta b&=&\beta_3+\alpha_3a. \nonumber
\end{eqnarray}

Summing these, we get $3\delta b=\beta_1+\beta_2+\beta_3+\delta a$. If $3b<a$, this is obviously impossible, implying that a binomial of this type can only occur when $3b\geq a$.

If $\alpha_1\geq 2$, then $\delta b-\beta_1\geq 2a>2\cdot(2b)=4b$. Therefore, $(\delta-4)b>\beta_1$, contradicting  $\delta\leq 3$ and $\beta_1\geq 0$. Thus, $\alpha_i\leq 1, i=1,2,3$.

If $\delta=3$, then $\alpha_i=1, i=1,2,3$, for otherwise, if $\alpha_1=0$ then $\beta_1=3b\geq a$. Thus, we retrieve $$H=w^3-(xyz)^{3b-a}tuv=E.$$
If $\delta=2$, then $\alpha_i=0,i=1,2,3$, for otherwise, if $\alpha_1=1$, then $\beta_1=2b-a<0$. Since $\delta=\alpha_1+\alpha_2+\alpha_3$ we obtain a contradiction.

\medskip

\noindent\underline{Type 2:} $H=x^{\beta_1}w^{\delta}-y^{\beta_2}z^{\beta_3}t^{\alpha_1}u^{\alpha_2}v^{\alpha_3}, \beta_1>0$. Then we have

\begin{eqnarray}
\beta_1+\delta b&=&\alpha_1a \nonumber \\
\delta b&=&\beta_2+\alpha_2a \nonumber\\
\delta b&=&\beta_3+\alpha_3a. \nonumber
\end{eqnarray}

If $\delta=2$, then $\alpha_2=\alpha_3=0$, otherwise $\beta_2<0$ (or $\beta_3<0$), as $a>2b$. Then $\alpha_1=2$, $\beta_1=2a-2b$, $\beta_2=\beta_3=2b$, giving $$H=x^{2a-2b}w^2-y^{2b}z^{2b}t^2=(x^{a-b}w-(yz)^bt)(x^{a-b}w+(yz)^bt) \in (L).$$

If $\delta=3$, then $\alpha_1\geq 1$. If $\alpha_1=1$, then $a>3b$, as otherwise $\beta_1\leq 0$. Then $\alpha_2=\alpha_3=0$, as otherwise $\beta_2=3b-\alpha_2a<0$ (or $\beta_3<0$). But this is still a contradiction since $\delta=\alpha_1+\alpha_2+\alpha_3$.

If $\delta=3$ and $\alpha_1=2$, then we can assume $\alpha_2=1$ (consequently, $3b>a$) and $\alpha_3=0$, leading to $$H=x^{2a-3b}w^3-y^{3b-a}z^{3b}t^2u=x^{2a-3b}(w^3-(xyz)^{3b-a}tuv)+z^{3b-a}y^{3b-a}tu(x^av-z^at).$$

If $\delta=3$ and $\alpha_1=3$, then $\alpha_2=\alpha_3=0$, and so $$H=x^{3a-3b}w^3-y^{3b}z^{3b}t^3=(x^{a-b}w-(yz)^bt)(x^{2a-2b}w^2+x^{a-b}(yz)^btw+(yz)^{2b}t^2)\in (L).$$

We conclude that there exist no binomials of this type that are not contained in $(L)$.

\medskip

\noindent\underline{Type 3:} $H=x^{\beta_1}y^{\beta_2}w^{\delta}-z^{\beta_3}t^{\alpha_1}u^{\alpha_2}v^{\alpha_3},\beta_1,\beta_2>0$. Then we have
\begin{eqnarray}
\beta_1+\delta b&=&\alpha_1a \nonumber \\
\beta_2+\delta b&=&\alpha_2a \nonumber\\
\delta b&=&\beta_3+\alpha_3a. \nonumber
\end{eqnarray}

As $\beta_1,\beta_2>0$, we must have $\alpha_1,\alpha_2\geq 1$.

If $\delta=2$, then $\alpha_1=\alpha_2=1$ and $\alpha_3=0$, thus retrieving
$$H=x^{a-2b}y^{a-2b}w^2-z^{2b}tu=H_1.$$

If $\delta=3$, then $\alpha_1=\alpha_2=\alpha_3=1$, and hence $\beta_1=a-3b>0$ and $\beta_3=3b-a\geq 0$; a contradiction. Or, $\alpha_1=2,\alpha_2=1,\alpha_3=0$, giving $\beta_1=2a-3b$ and $\beta_2=a-3b>0$. So $\beta_1>a$, a situation that can be disregarded from the beginning.

\medskip

\noindent\underline{Type 4:} $H=x^{\beta_1}y^{\beta_2}z^{\beta_3}w^{\delta}-t^{\alpha_1}u^{\alpha_2}v^{\alpha_3}, \beta_1,\beta_2,\beta_3>0$. Then we have
\begin{eqnarray}
\beta_1+\delta b&=&\alpha_1a \nonumber \\
\beta_2+\delta b&=&\alpha_2a \nonumber\\
\beta_3+\delta b&=&\alpha_3a. \nonumber
\end{eqnarray}

This gives $\alpha_i\geq 1, i=1,2,3$, and since $\delta\leq 3$, one must have $\alpha_i=1,i=1,2,3$. This means $\beta_i=a-3b, i=1,2,3$, and $a>3b$. Therefore, we retrieve $$H=(xyz)^{a-3b}w^3-tuv=E'.$$

\bigskip

(b) The argument to prove almost Cohen--Macaulayness will proceed via mapping cones, in a way reminiscent of the proof of Theorem~\ref{main2}.

First we compute the colon ideals between iterated subideals of $\mathcal I\subset S=R[t,u,v,w]$ in a suitable order as to profit from the nature of the Sylvester generators in part (a).

For lighter reading below we will sometimes write $L$ instead of $(L)$.

We fix once for all the monomial lexicographic order with $w>t>u>v>x>y>z$ to be used in the line of argument to follow.

\medskip

{\sc Claim 1:} $(L):H_1=(x^b,y^b,z^{a-b})S$.

Since $z^{a-b}H_1=(xy)^{a-2b}wg_3+y^{a-b}vg_1+z^btf_3$, and since $(g_1,g_2):H_1=(x^b,y^b)S$, the inclusion $\supseteq$ is obvious. For the reverse inclusion, under the above monomial order it is easy to see that $L$ is a Gr\"{o}bner basis, and therefore the initial ideal of $(L)$ is $${\rm in}(L)=(x^{a-b}w,y^{a-b}w,z^{a-b}w,y^at,z^at,z^au).$$
On the other hand, ${\rm in}(H_1)=(xy)^{a-2b}w^2$. It easily follows that ${\rm in}(L):{\rm in}(H_1)=(x^b,y^b,z^{a-b})S$.
But since the right side in the colon $(L):H_1$ is a principal ideal, we have ${\rm in}(L:H_1)\subset {\rm in}(L):{\rm in}(H_1)$.
Therefore, ${\rm in}(L:H_1)\subset (x^b,y^b,z^{a-b})S$.

We repeat this scheme for $f-{\rm in}(f)$, and so forth, resulting that every monomial in the expression of an arbitrary $f\in (L):H_1$ is in fact in the ideal $(x^b,y^b,z^{a-b})S$.

\medskip

{\sc Claim 2:} $(L,H_1):H_2=(x^b,y^{a-2b},z^b)S$.

The inclusion $\supseteq$ is clear once we observe that $y^{a-2b}H_2=z^{a-2b}H_1-tf_3$. By a similar token as above it obtains
$${\rm in}(L,H_1)=(x^{a-b}w,y^{a-b}w,z^{a-b}w,y^at,z^at,z^au,(xy)^{a-2b}w^2)S
,$$
and hence ${\rm in}(L,H_1):\underbrace{(xz)^{a-2b}w^2}_{{\rm in}(H_2)}=(x^b,y^{a-2b},z^b)S$. We conculde as in the first claim.

\medskip

{\sc Claim 3:} $(L,H_1,H_2):H_3=(x^{a-2b},y^b,z^b)S$.
The proof is a repeat of the scheme of the proofs of the previous two claims.

\medskip

{\sc Claim 4:} (Case $3b>a$) $(L,H_1,H_2,H_3):E=(x^{a-b},y^{a-b},z^{a-b},(xy)^{a-2b},(xz)^{a-2b},(yz)^{a-2b})S$.

Same scheme of proof with the required changes.
Moreover, the case $a\geq 3b$ is handled in the same manner, with the obvious modifications; note also that $(L,H_1,H_2,H_3):E'=(x^b,y^b,z^b)^2$.

\medskip

We now proceed to the construction of the successive mapping cones, starting from the free resolution of $S/(L)$.
Since the ideal $I\subset R$ is an almost complete intersection of finite length, $S/(L)$ is Cohen--Macaulay (\cite[Corollary 10.2]{Trento}). Its codimension is at least $3$ (the codimension of $S/{\mathcal I}$), but since $L\subset (x,y,z)S$ then the codimension is $3$.

Let $0 \rar S^{\beta_3} \lar  S^{\beta_2} \lar S^6 \lar S$ stand for a free resolution of $S/L$ (actually, one can take $\beta_2=7$ and, consequently, $\beta_3=2$, but we don't need this extra information).

Let $0\rar S\lar S^3 \lar S^3 \lar S$ denote the minimal free resolution of $(L):H_1=(x^b,y^b,z^{a-b})S$ extended from $R$ by flat base change $R\subset S$.
Consider the map of complexes induced by multiplication by $H_1$ on $S$:

$$\begin{array}{ccccccccccc}
0 & \rar & S^{\beta_3}  & \lar & S^{\beta_2} & \lar & S^6 & \lar & S & \rar & 0\\
&&\uparrow && \uparrow && \uparrow && \uparrow && \\
0 & \rar & S  & \lar & S^3 & \lar & S^3 & \lar & S & \rar & 0
\end{array}.
$$
Then its mapping cone
\begin{equation}\label{mcone1}
0 \rar S\lar  S^{\beta_3+3} \lar  S^{\beta_2+3} \lar S^7 \lar S
\end{equation}
is a resolution of $S/(L,H_1)$ (it won't be minimal as there is a cancellation of a summand $S$, but again we don't need this additional information).

Proceed to the next step, by taking the map of complexes from the minimal free resolution of $(L,H_1):H_2=(x^b,y^{a-2b},z^b)S$ to (\ref{mcone1}) induced by  by multiplication by $H_2$ on $S$.
Since (\ref{mcone1}) has length $4$ and  $S/(x^b,y^{a-2b},z^b)S$ has homological dimension $3$, the resulting mapping cone resolves $S/(L,H_1,H_2)$ and is again of length $4$.

Likewise, iterating next with $H_3$, we arrive at a free resolution of the subideal $(L,H_1,H_2,H_3)$ of length $4$.

In order to get to the final step, we need to know a minimal free resolution of $S/Q$, where $Q=(L,H_1,H_2,H_3):E$ , if $3b>a$ (or $Q=(L,H_1,H_2,H_3):E'$, if $a\geq 3b$).
But again, since $Q$ is the extension to $S$ of an $(x,y,z)$-primary ideal of $R$, then $S/Q$ is Cohen--Macaulay of codimension $3$.

Therefore,  we are in the same situation as in the previous steps and hence $\Ree I=S/{\mathcal I}$ has a free resolution of length at most $4$.
This means that $\Ree I$ is almost Cohen--Macaulay.
\qed

\begin{Remark}\rm
(a) When $(R,\fm)$ is a regular local ring and $I$ is an $\fm$-primary ideal whcih is an almost complete intersection, there is an a priori precise relation between the reduction number with respect to reduction generated by a regular sequence and the relation type of $I$ (see \cite{MuPl}). Since we are dealing with forms of different degrees, we felt safer to give an independent self-contained argument
that does not assume that the relation type is $\leq 3$.

(b) The procedure opens a clear path towards the case of $n$ variables, as it is a lot more straightforward to guess the generators in the light of Sylvester forms and most likely get all essential generators of the Rees ideal.
Furthermore, iterating Sylvester forms leads to a stepwise procedure to reach a bound of the homological dimension of $\Ree I$.
Still, it can be a mighty task to organize the proof accordingly.
\end{Remark}



\begin{thebibliography}{99}

\bibitem{Dandrea}{T. C. Benitez and C.  D'Andrea, The Rees Algebra of a monomial plane parametrization, arXiv:1311.5488v2 [math.AC].}

\bibitem{CHW}{D. Cox, J. W. Hoffman and H. Wang, Syzygies and the Rees
algebra, J. Pure Appl. Algebra, {\bf 212} (2008), 1787--1796.}


\bibitem{aha}{A. V. D\'oria, S. H. Hassanzadeh and A. Simis, A characteristic free criterion of birationality, Advances in Math., {\bf 230} (2012), 390--413.}


\bibitem{HS}{S. H. Hassanzadeh and A. Simis, Plane Cremona maps:
saturation and regularity of the base ideal, J. Algebra, {\bf 371} (2012), 620–-652.}

\bibitem{Trento}{J. Herzog, A. Simis and W. V. Vasconcelos, Koszul homology and blowing-up rings, in {\sc Commutative Algebra} (S. Greco and G. Valla, eds.),  Lecture Notes in Pure and Applied Math.  {\bf 84}, Marcel-Dekker, New York, 1983, 79--169.}

\bibitem{syl1}{J. Hong, A. Simis, W. V. Vasconcelos, On the homology of two-dimensional elimination, J. Symb. Comp. {\bf 43} (2008) 275--292.}

\bibitem{syl2}{J. Hong, A. Simis, W. V. Vasconcelos, The equations of almost complete intersections, Bull. Braz. Math. Soc. (New Series), {\bf 43} (2012), 171--199.}


\bibitem{syl3}{J. Hong, A. Simis, W. V. Vasconcelos,  Extremal Rees algebras, J. Comm. Algebra, {\bf 5} (2013), 231--267.}

\bibitem{HM}{S. Huckaba and T. Marley, Hilbert coefficients and the depths of associated graded rings, J. London Math. Soc. {\bf 56} (1997), 64--76.}

\bibitem{MuPl}{F. Mui\~nos and F. Planas-Vilanova, Equations of powers of equimultiple ideals of deviation one, Proc. Amer. Math. Soc. {\bf 141} (2013), 1241--1254.}

\bibitem{RoSw}{M. E. Rossi and I. Swanson, Notes on the behavior of the Ratliff-Rush filtration, {\em in} Contemp. Math., 331, Amer. Math. Soc., Providence, RI, 2003.}


\bibitem{Vas}{W. Vasconcelos, {\sc Computational Methods in Commutative Algebra and Algebraic Geometry}, Springer-Verlag, New York 2004.}



\end{thebibliography}
\end{document}